\documentclass{amsart}
\newtheorem{Theorem}{Theorem}[section]
\newtheorem{Lemma}[Theorem]{Lemma}
\newtheorem{Proposition}[Theorem]{Proposition}
\newtheorem{Corollary}[Theorem]{Corollary}

\theoremstyle{definition}
\newtheorem{Definition}[Theorem]{Definition}

\newtheorem{ex}[Theorem]{Example}

\theoremstyle{remark}
\newtheorem{remark}[Theorem]{Remark}

\theoremstyle{remark}

\theoremstyle{remark}

\theoremstyle{remark}
\newtheorem{case}[Theorem]{Case}

\theoremstyle{remark}

\numberwithin{equation}{section}

\newcommand{\cal}{\mathcal}
\renewcommand{\O}{{\mathcal O}}
\renewcommand{\o}{{\mathcal O}}

\newcommand{\I}{{\mathcal I}}
\newcommand{\E}{{\mathcal E}}
\renewcommand{\L}{{\mathcal L}}
\newcommand{\Proj}{{\mathbb P}}

\newcommand{\R}{{\mathbb R}}
\newcommand{\HH}{{\mathbb H}}
\newcommand{\OO}{{\mathbb O}}
\newcommand{\G}{{\mathbb G}}

\newcommand{\Z}{{\mathbb Z}}

\newcommand{\C}{{\mathbb C}}
\newcommand{\p}{{\mathbb P}}
\renewcommand{\a}{{\`a}}
\newcommand{\f}{\varphi}

\newcommand{\map}{\dasharrow}

\newcommand{\ra}{\rightarrow}

\newcommand{\flip}{{\stackrel{\scriptstyle - - >}{  }}}

\begin{document}

\title{Varieties with one apparent double point}

%author one
\author{Ciro Ciliberto}
\address{Dipartimento di Matematica, Universit\'a di Roma Tor Vergata,
Via della
Ricerca Scientifica,
 00133 Roma, Italia}
\email{cilibert@axp.mat.uniroma2.it}

%author two
\author{Massimiliano Mella}
\address{Dip. di
Matematica, Universit\a di Ferrara, Via Machiavelli 35,
 44100 Ferrara, Italia}
\email{mll@unife.it}

%author three
\author{Francesco Russo}
\address{Departamento de Matematica, Universidade Federal
de Pernambuco, Cidade Universitaria, 50670-901 Recife--PE,
Brasil}
\email{frusso@dmat.ufpe.br}

\begin{abstract}
 The number of apparent double points of a smooth, irreducible
projective variety $X$ of dimension $n$ in $\Proj^{2n+1}$ is the
number of secant lines to $X$ passing through the general point of
$\Proj^{2n+1}$. This classical notion dates back to Severi. In
the present  paper we classify smooth varieties of dimension at
most three having one apparent double point. The techniques
developed for this purpose allow  to treat a wider class of
projective varieties.
\end{abstract}

\maketitle

\section*{Introduction}

The number $\nu(X)$ of {\it apparent double points} of
 a smooth, irreducible projective variety $X$ of dimension $n$ in
$\Proj^{2n+1}$ is the number of secant lines to $X$ passing
through the general point of $\Proj^{2n+1}$. The variety $X$ is
called a {\it variety with one apparent double point}, or {\it
OADP}-variety, if $\nu(X)=1$. Hence {\it OADP}-varieties of
dimension $n$ are not defective, and actually they can be
regarded as {\it the simplest non defective varieties} of
dimension $n$ in $\Proj^{2n+1}$. Probably this is the reason why
{\it OADP}-varieties attracted since a long time the attention of
algebraic geometers. There are other reasons however. Among which
we will mention here the fact that the known {\it OADP}-varieties
are rather interesting objects, as we will indicate later on in
this paper, and the unexpected relation of {\it OADP}-varieties
with an important subject like Cremona transformations (see \cite
{AR}). Unpublished work of F. Zak also suggests connections with
the classification of self-dual varieties. \par

The first instance of an {\it OADP}-variety is the one of a rational
normal cubic in $\Proj^3$, which actually
is the only {\it OADP}-curve. It goes back to Severi \cite {Se} the
first
attempt of classifying {\it OADP}-surfaces. According to him the only
{\it OADP}-surfaces in $\Proj^5$ are
quartic rational normal scrolls and del Pezzo quintics (see Theorem \ref
{OAPDsurfaces} below). Severi's
nice geometric argument contained, however, a gap, recently fixed by
the third author of the present
paper \cite {Ru}. Meanwhile, inspired by Severi's result, Edge produced
in \cite {Edge} two infinite series of
{\it OADP}-varieties of any dimension (see Example \ref {edgevarieties}). In
dimension $2$ the two types of {\it
OADP}-surfaces are either Edge varieties or deformations of them (see again
Example \ref {edgevarieties}). In more recent times
E. Sernesi called the attention on Edge's
paper, asking  whether, like in the surface case, all {\it OADP}-varieties can be
somehow recovered by Edge varieties. F. Zak has been the first who produced, in
dimension $4$ or more a few very interesting sporadic examples of
{\it OADP}-varieties (see \S \ref {examples} below). Unfortunately
however, there is no other infinite list of
examples but Edge's varieties and rational normal scrolls. \par

The challenge of classifying {\it OADP}-varieties seems to be quite
hard. A first remark is that one may
associate to an {\it OADP}-variety $X$ of dimension $n$ another variety
$F(X)$ in $\Proj^{2n+1}$ which is
called the {\it focal variety} of $X$. Basically $F(X)$ is the locus of
points in $\Proj^{2n+1}$ which are
contained in more than one, hence in infinitely many, secant lines to
$X$ (see \S \ref {definitions} below). Of
course $X$ is contained in $F(X)$ and one proves that $F(X)$ has at most
dimension $2n-1$. One is therefore
tempted to classify {\it OADP}-varieties according to the dimension of
their focal loci. This essentially works
in the case of surfaces, but falls short for higher dimensions, where
one needs some new idea.\par

In the present paper we deal with the case of {\it OADP}-threefolds, and
we are able to give a full
classification theorem for them (see Theorem \ref {Th:classification}).
The result is that, like in the case of
surfaces, the degree $d$ of such varieties is bounded and precisely one
has $5\leq d\leq 8$. The case $d=5$
corresponds to rational normal scrolls, the cases $d=6,7$ correspond to
Edge's varieties, whereas the case $d=8$
corresponds to an interesting scroll in lines with rational hyperplane
section and sectional genus $3$, which
naturally arises in the context of Cremona transformations (see \cite
{AR}, \cite {HKS}). \par

The new idea on which we base our analysis, and which we believe could
be useful also in other contexts, is to
study the degeneration of the singular locus of a general projection of
a variety when the center of projection
tends to meet the variety itself. Using this, we are able to prove a few
interesting general facts about {\it
OADP}-varieties. What is really relevant for us here, is that any
irreducible $1$-dimensional
component of the intersection of an {\it OADP}-variety of dimension $n$
with $n-1$ general tangent hyperplanes at
a general point is a rational curve. Applying this to surfaces one
recovers Severi's classification theorem right
away. In the case of threefolds the situation is much more involved and
requires a case by case rather delicate
analysis. Nonetheless it works, and at the end it provides the bound
$d\leq 9$ for the degree. Once one has this,
the conclusion of the classification can be done in several ways: we
used in a rather efficient way adjunction
theory. This enables us to exclude the case $d=9$ and to conclude that
the other cases $5\leq d\leq 8$
correspond to the aforementioned varieties.\par

It should be noticed that M. G. Violo, in her Ph.D. thesis \cite
{violo}, also studied the same classification
problem. She tried to adapt to higher dimensions Severi's argument for
surfaces. However, in order to do so,  she
had to make a rather unnatural assumption and her list of {\it
OADP}-threefolds does not include the degree $8$
scroll. \par

Our approach could in principle be applied also in higher dimensions.
However the analysis becomes then
much heavier, to the extent of being in fact discouraging. What it
suggests is that it should be
possible to bound the degree of an {\it OADP}-variety of dimension $n$
by a function depending on $n$. In other
words, it should still be the case that, in any dimension, one has
finitely many families of {\it
OADP}-varieties. The bound that our approach suggests for the degree is
exponential in $n$, and this is
one of the reasons why the classification, in this way, becomes rather
unfeasible. It could happen,
however, that the cases which are really possible are much more
restricted: after all, even in dimension $3$, the
case $d=9$, a priori possible according to our approach, does not in
fact exist.\par

One more information which we are able to give is that both an {\it
OADP}-variety and its general hyperplane
section  have Kodaira dimension $-\infty$ (see Proposition \ref {RPT}).
This suggests an alternative approach to the
classification which we briefly outline in Remark \ref
{rationalhypsect}. It also relates to an interesting
conjecture of Bronowski \cite {Br} to the effect that a variety $X\subset
\Proj^{2n+1}$ of dimension $n$ has one apparent double point if and only if the
projection of $X$ to $\Proj^n$ from the projective tangent space to $X$ at its general point
is birational  (see Remark \ref
{Bronowski}). Here we prove that, in general, the degree of aforementioned
projection is bounded above by the number of apparent double points of $X$, which
yields one implication of Bronowski's conjecture (see Proposition \ref
{prebronowski} and Corollary \ref {bronowskiconj}). In particular we have the
important consequence that {\it OADP}-varieties are
rational.   \par

Finally, as an application of our classification of {\it
OADP}-threefolds, we are able to classify  {\it
OADP}-varieties which are  Mukai varieties, i.e. such that their
general curve
sections are canonical curves (see Theorem \ref {Th:Mukai}). Related
classification results for {\it
OADP}-varieties of {\it small degree} can be found in \cite {AR}).\par

The paper is organized as follows. In section \ref {definitions} we give
the main definitions and prove a few
basic facts about {\it OADP}-varieties. In \S \ref{examples} we present
all the known examples of {\it
OADP}-varieties. In \S \ref {degenerations} we make our analysis of
degenerations of projections. In \S \ref
{applications} we apply it to the study of a special class of varieties,
which we call {\it varieties
of full rational type}, which include the {\it OADP}-varieties. As a
first result we give a very easy proof of
Severi's classification theorem of {\it OADP}-surfaces. In \S \ref
{tanbeh} we use the previous analysis for the
study of the tangential behavior of {\it OADP}-threefolds, i.e. we
study the intersection of such a threefold with
its general tangent space. As a result we are able to prove the bound
$5\leq d\leq 9$ for the degree $d$ of a {\it
OADP}-threefold.  On the way, we need a few technical propositions, some
of them of independent interest, which we
collected in \S \ref {lemmata}. In \S \ref  {classification} we prove
our classification theorem for {\it
OADP}-threefolds. Finally in \S \ref {mukai}
we prove the classification of {\it OADP}-Mukai varieties. \medskip

\noindent {\bf Acknowledgments:} The authors wish to thank L. Chiantini,
E. Sernesi and F. Zak for several
interesting discussions on the subject of this paper. The first and
second author have been supported by the EC
research project HPRN-CT-2000-00099, EAGER. The third author has been
supported by the CNPq grant 300761/97 and PRONEX-ALGA (Algebra Comutativa e Geometria Algebrica).\par

\section{Preliminaries}\label {definitions}

Let $X\subset\Proj^r$ be an irreducible, non-degenerate, projective
variety
of dimension $n\geq 1$ over ${\C}$. If $d$ is the degree of $X$ we will
sometimes denote $X$ by $X_{n,d}$ or
by $X_n$.\par

Let $X(2)$ be the twofold symmetric
product of $X$. One considers ${\cal S}(X)$ the {\it abstract secant
variety} of $X$, i.e. the Zariski closure of
the scheme

$$I(X):=\{([x,y],p)\in X(2)\times \Proj^{r}: x\neq y, p\in <x,y>\}$$

\noindent inside $X(2)\times \Proj^{r}$. The image of the projection of
${\cal S}(X)$ to $\Proj^{r}$ is the {\it
secant variety} $Sec(X)$ of $X$. Thus one has the projection map $p_X:
{\cal S}(X)\to Sec(X)$. Since $\dim({\cal S}(X))=2n+1$, then
$\dim(Sec(X))\leq \min \{2n+1,r\}$, and the variety $X$ is said to be
{\it defective} if the strict inequality
holds. If $x,y\in X$ are distinct smooth points, the  {\it secant line}
$L=<x,y>$ is contained in $Sec(X)$. As
limit of secant lines, one has also the {\it tangent lines}, i.e. the
lines contained in the tangent space
$T_xX$ and passing through $x$, where $x\in X$ is a smooth point. The
Zariski closure of the union of the tangent
lines is the so-called {\it tangential variety} $Tan(X)$ of $X$, which
is contained in $Sec(X)$.  The tangent lines
are considered to be {\it improper} secant lines, whereas the usual
secant lines are sometimes called {\it
proper} secants. \par

Let us assume now that $X$ is smooth and that $r\geq 2n+1$. Let
$\epsilon$ be
an integer such that $0\leq \epsilon\leq n-1$ and let $\Pi$ be a general
$\Proj^s$
in $\Proj^r$, with $s=r-n-\epsilon-2$. We denote by $\pi_\epsilon: X\to
\Proj^{n+1+\epsilon}$ the
projection of $X$ from $\Pi$, by $X^{(\epsilon)}$ its image and by
$\Gamma^{(\epsilon)}$ the singular scheme of $X^{(\epsilon)}$. We will
assume from now on
that $X$ is not defective, i.e. that its secant variety $Sec(X)$ has
dimension $2n+1$. Then
$\Gamma^{(\epsilon)}$ has pure dimension $n-\epsilon-1$.\par

The projective invariants of the varieties $\Gamma^{(\epsilon)}$, which
are projective invariants of the
variety $X$ itself, are relevant in various classification problems. For
instance, if $\epsilon=n-1$, then
$\Gamma^{(n-1)}$ consists of finitely many points, which are {\it
improper double points} of $X^{(n-1)}\subset
\Proj^{2n}$, i.e. double points whose tangent cone consists of two
subspaces of dimension $n$ meeting transversely
at the points in question. Following a classical terminology, we call
the points of $\Gamma^{(n-1)}$ the {\it
apparent double points} of $X$. Their number, which is not $0$ since $X$
is not defective, will be denoted by
$\nu(X)$. Usually one sets $\nu(X)=0$ for defective varieties $X$. \par

It is useful to point out  right away the following basic facts:

\begin{Proposition} \label {linorm} Let $X\subset \Proj^r$ be a smooth,
irreducible, non-degenerate,
non-defective, projective variety of dimension $n$, with $r\geq 2n+1$.
One has:
\begin{itemize}
\item[i)] if $r=2n+1$ then $\nu(X)$ is the number of secant lines to
$X$ passing through a
general point of $ \Proj^{2n+1}$,
\item[ii)] if $r>2n+1$ then $\nu(X)=\deg(Sec(X))\cdot\deg(p_X)$ and
$r\leq
2n+\nu(X)-1$, hence $\nu(X)\geq 3$,
\item[iii)] if $\nu(X)\leq 2$ then $r=2n+1$ and $X$ is linearly
normal.
\end{itemize}
\end{Proposition}

\begin{proof}  Assertion $(i)$ is trivial, and $(iii)$ is an immediate
consequence of $(ii)$, which we are
now going to prove. The first assertion of $(ii)$ is also trivial and
one has:

$$\nu(X)=\deg(Sec(X))\cdot \deg(p_X)\geq \deg(Sec(X))\geq
r-\dim(Sec(X))+1=r-2n$$

Assume the equality holds. Then $Sec(X)$ would be a variety of minimal
degree (see \cite{EH} and Example \ref{scrolls} below). Such a variety, if singular, is a cone. Actually
$Sec(X)$ is singular along $X$, hence $Sec(X)$
would be a cone, whose vertex would contain the span of $X$, which is
the whole $\Proj^r$, a contradiction.
\end{proof}

A smooth, irreducible, non-degenerate,
non-defective, projective variety $X\subset \Proj^r$ of dimension $n\geq
1$ is called a variety with {\it
one apparent double point}, or a {\it OADP-variety}, if $\nu(X)=1$. Note
that, by part $(iii)$ of Proposition \ref{linorm}, if
$X$ is  an  {\it OADP}-variety, one has $r=2n+1$ and $X$ is linearly
normal in $\Proj^{2n+1}$ (see also
\cite {Se} and \cite{Ru}). \par

To say that $X$ is an {\it OADP-variety} is the same as saying that the
projection $p_X: {\cal S}(X)\to
\Proj^{2n+1}$ is dominant and birational. This fact has some important
consequences.\par

Let us consider, in general, $X_n\subset \Proj^{2n+1}$ an irreducible,
non-degenerate, non-defective, projective
variety. We can see ${\cal S}(X)$ as a family of dimension $2n$ of lines
in $\Proj^{2n+1}$. Let $L=<x,y>$, $x,y\in
X$, be a proper secant line. In this situation there is {\it focal}
square matrix $F_L$ arising, whose
rows are given by the sections of the normal bundle $N_{L,{\Proj^{2n+1}}}$
corresponding to the $2n$ independent
infinitesimal deformations of $L$ determined by $2n$ independent vectors
in $T_{[x,y]}(X(2))$ (see \cite {cs} for
the general theory of foci, which we are going to apply in this paper).
Since $X$ is not defective, hence $p_X:
{\cal S}(X)\to \Proj^{2n+1}$ is dominant, then $\det(F_L)$ is not
identically zero for $L$ general. The scheme
$Z_L$ of degree $2n$ defined by the equation $\det (F_L)=0$, if
$\det(F_L)$ is not identically zero, coincides with
$nx+ny$. The same considerations apply if $L$ is an improper secant
line. In this case $L$ is tangent to $X$ at a
smooth point $x\in X$ and the focal matrix $F_L$ can still be considered.
If $\det(F_L)$ is not identically zero,
the focal scheme $Z_L$ of degree $2n$ defined by $\det (F_L)=0$ is
$2nx$. A proper or improper secant line $L$ is
called a {\it focal line} if $\det (F_L)\equiv 0$. The union of focal
lines is a Zariski closed subset $F(X)$ of
$\Proj^{2n+1}$ called the {\it focal locus} of $X$ and its points are
called {\it foci}. \par

Assume now $X$ is an {\it OADP}-variety. By the general theory of foci,
$F(X)$ coincides with the
indeterminacy locus of the inverse map $ p_X^{-1}: \Proj^{2n+1}\dasharrow {\cal
S}(X)$. Hence:
\begin{itemize}

\item[i)] $z$ is a focus if and only if there are two distinct,
hence infinitely many, secant lines through
$z$,
\item[ii)]the secant line $L$ is a focal line if and only if for
some $z\in L-L\cap X$ there is
some other secant line, different from $L$, containing $z$. This happens
if and only if for  every $z\in L$, there
are infinitely many secant lines containing $z$.
\end{itemize}

\begin{remark} \label{singularOADP} The definition of {\it
OADP}-varieties can be extended also to singular
varieties. A non-degenerate, projective variety $X$ of pure dimension
$n$ in $\Proj^{2n+1}$ has one
apparent double point if there is a unique secant to $X$ passing through
the general point of $\Proj^{2n+1}$.
These varieties share with their smooth relatives, several properties.
For example they also turn out to be
linearly normal (see the proof of Proposition \ref {linorm}). Though the
subject of singular {\it OADP}- varieties
is very interesting, we will not consider it in this paper.\end{remark}
\medskip

Going back to the loci $\Gamma^{(\epsilon)}$ for general varieties, let
us consider the next one, corresponding to
the case $\epsilon=n-2$, thus we have to assume $n\geq 2$. Then
$\Gamma^{(n-2)}$ is a curve whose geometric genus
we denote by $PG(X)$ and we will call it the {\it projective genus} of
$X$. Then
$PG(X)\geq 0$, and we will say that $X$ is of {\it rational projective
type}, or a
{\it RPT-variety}, if
$PG(X)=0$ and, in addition, the curve $\Gamma^{(n-2)}$ is irreducible.
More generally, we will say that
$X$ is of {\it rational projective $i$-type}, for some $i=1,...,n-1$, or
a {\it
RPT$_i$-variety}, if the variety $\Gamma^{(n-i-1)}$ is an irreducible,
rational,
$i$-dimensional variety. We will say that $X$ is of {\it full rational
projective
type}, or a {\it FRPT-variety}, if it is of {\it rational projective
$i$-type}, for
all $i=1,...,n-1$. \par

As the following proposition shows, examples of $FRPT$-varieties are
provided by {\it OADP-varieties}.

\begin{Proposition} \label {OADP-FRPT} Every {\it OADP}-variety is a
{\it
FRPT}-variety.
\end{Proposition}

\begin{proof}  For every $0\leq \epsilon\leq n-1$ we consider $\Pi$ a
general
$\Proj^{n-\epsilon-1}$ and $\pi_\epsilon: X\to \Proj^{n+1+\epsilon}$ the
projection
of $X$ from $\Pi$. Since through the general point $p\in \Pi$ there is a
unique
secant line to $X$ corresponding to a point of
$\Gamma^{(\epsilon)}$, we have a rational map $f_\epsilon: \Pi\to
\Gamma^{(\epsilon)}$,
which is birational. \end{proof}

\begin{remark} \label {rationality} If $X\subset \Proj^{2n+1}$ is a {\it
OADP}-variety, then its
symmetric product $X(2)$ is rational. Indeed, if $H$ is a general
hyperplane of $\Proj^{2n+1}$, then, referring
to the map $\pi: {\cal S}(X)\to
\Proj^{2n+1}$, one has that $H'=\pi^*(H)$ is rational and the projection
map $H'\to X(2)$ is clearly birational.
Thus the Kodaira dimension of $X$ is $\kappa(X)=-\infty$ and
$h^i(X,{\cal O}_X)=0$, $1\leq i\leq n$. This can be proved as in
\cite{mattuck}.
\label{mattuck}
\end{remark}
\medskip

We prove now a lemma indicated to us by F. Zak, which will be useful
later on. First we recall a
definition. Let $X\subset \Proj^r$ be an irreducible, projective variety
of dimension $n$. For a point
$p\in (Sec(X)\setminus X)$ we consider the Zariski
closure $C_p(X)$ of the union of (proper or
improper) secant lines to $X$ passing through $p$. We denote by $S_p(X)$
the scheme theoretical intersection of
$C_p(X)$ and $X$, and we call $S_p(X)$ the {\it entry locus} of $p$ with respect to
$X$. Intuitively, this is the Zariski closure of the set of all points $x\in X$
for which there is a point $y\in X$ such that $p\in\langle x,y\rangle$.\par

Assume $Sec(X)=\Proj^r$,
which yields $r\leq 2n+1$. A count of parameters shows that for $p\in
\Proj^r$ a general point, $S_p(X)$
has  dimension $2n+1-r$. \par

Notice that $X_n\subset \Proj^{2n+1}$ is an {\it OADP}-variety
if and only if the entry locus $S_p(X)$ of a general point $p\in
\Proj^{2n+1}$ is a pair of points.  \par

\begin{Lemma} \label{zak} Let $Y\subset \Proj^N$ be an irreducible,
non-degenerate, projective
variety of dimension $m$, such that $Sec(Y)=\Proj^N$. Set $n=N-m-1$ and
$r=2n+1$. Consider $\Pi$ a $\Proj^r$ in
$\Proj^N$ intersecting $Y$ along a smooth variety $X$ of dimension $n$.
Suppose that for $p\in \Pi$ general the
entry locus $S_p(Y)$ is a quadric of dimension $2m+1-N$. Then $X$ is an
{\it OADP}-variety.\end{Lemma}

\begin{proof} Let $p$ be a general point of $\Pi$. The entry locus
$S_p(X)$ is the intersection of the quadric
$S_p(Y)$ with  $\Pi$. By a dimension count, this intersection is not
empty.  This proves that
$Sec(X)=\Pi$ and therefore $S_p(X)$ is finite. Since, as we saw, it is
the intersection of the quadric $S_p(Y)$
with $\Pi$, it consists of two points. This implies the assertion about
$X$.\end{proof}

To put the previous lemma in perspective the reader should perhaps take
into account the relations
between {\it OADP}-varieties and Cremona transformations, as indicated in
\cite {AR}. We will not dwell on this
here but we will only mention the following result (see \cite {AR}):

\begin{Proposition}\label{cremona} Let $X\subset \Proj^{2n+1}$ be a
smooth, irreducible, non-degenerate,
projective variety, which is not defective. Suppose that the linear
system of quadrics
$|{\cal I}_{X,\Proj^{2n+1}}(2)|$, restricted to a general hyperplane of
$\Proj^{2n+1}$ defines a rational map
which is birational onto its image. Then $X$ is an {\it
OADP}-variety.\end{Proposition}

\begin{proof} Let $\phi$ be the map associated to $|{\cal
I}_{X,\Proj^{2n+1}}(2)|$. The map $\phi$ contracts any
secant line to $X$ to a point. Hence the general fiber of $\phi$ has
positive dimension $l$ and degree $d$. Since
the restriction of $\phi$ to a general hyperplane has to be birational
onto its image, we must have
$l=d=1$. This implies the assertion.\end{proof}

\section {Examples of {\it OADP}-varieties}\label {examples}

In this section we collect some examples, essentially the only ones
known to us, of smooth {\it
OADP}-varieties. We refer the reader to \cite {AR} for more details.\par

\begin{ex}\label {scrolls} ({\it Scrolls}). Let $0\leq a_0\leq a_1\leq
...\leq a_k$ be integers and set
$r=a_0+...+a_k+k$. Recall that a {\it rational normal scroll}
$S(a_0,...,a_k)$ in $\Proj^r$ is the image of the
projective bundle $\Proj:=\Proj(a_0,...,a_k):=\Proj({\cal
O}_{\Proj^1}(a_0)\oplus...\oplus{\cal
O}_{\Proj^1}(a_k))$ via the linear system $|{\cal O}_{\Proj}(1)|$. The
dimension of $S(a_0,...,a_k)$ is $k+1$,
its degree is $a_0+...+a_k=r-k$ and $S(a_0,...,a_k)$ is smooth if and
only if $a_0>0$. Otherwise, if
$0=a_0=...=a_i<a_{i+1}$, it is  the cone over $S(a_{i+1},...,a_k)$ with
vertex a $\Proj^i$.

One uses the simplified
notation $S(a_0^{h_0},...,a_m^{h_m})$ if $a_i$ is repeated $h_i$ times,
$i=1,...,m$.\par

Recall that rational normal scrolls, the (cones over) Veronese surface in
$\Proj^5$, and quadrics,  can be characterized as those non degenerate irreducible
varieties in a projective space having minimal degree (see \cite{EH}).\par

Given positive integers $0<m_1\leq ...\le m_h$ we will denote by
$Seg(m_1,...,m_h)$ the Segre embedding of
$\Proj^{m_1}\times ...\times \Proj^{m_h}$ in $\Proj^N$,
$N=(m_1+1)...(m_h+1)-1$. We use the shorter notation
$Seg(m_1^{k_1},...,m_s^{k_s})$ if $m_i$ is repeated $k_i$ times,
$i=1,...,s$. Recall that
$Pic(Seg(m_1,...,m_h))\simeq \Z^h$, generated by the line bundles
$\xi_i=pr_i^*({\cal O}_{\Proj^{m_i}}(1))$,
$i=1,...,h$. A divisor $D$ on  $Seg(m_1,...,m_h)$ is said to be of ${\it
type} (\ell_1,...,\ell_h)$ if ${\cal
O}_{Seg(m_1,...,m_h)}(D)\simeq \xi_1^{\ell_1}\otimes ... \otimes
\xi_h^{\ell_h}$. The hyperplane divisor of
$Seg(m_1,...,m_h)$ is of type $(1,...,1)$.\par

Notice now that $S(1^n)\subset {\Proj^{2n-1}}$ coincides with
$Seg(1,n-1)$. Let us point out the following easy
fact.

\begin{Lemma} \label {basiscrolls} Let $X\subset\Proj^{2n+1}$ be a
smooth, linearly normal,
regular variety of dimension $n$. Then $X$ is a scroll over a curve
if and only if $X$ is isomorphic either to
$S(1^{n-2},2^2)$, $n\geq 2$, or
$S(1^{n-1},3)$, $n\geq 1$. Moreover
for $n\geq 3$ these scrolls can be realized as divisors of type
$(2,1)$ on $Seg(1,n)$, the general one being
$S(1^{n-2},2^2)$, while the scroll $S(2^2)$ can be also realized as a
divisor of type $(0,2)$ on $Seg(1,2)$.
Furthermore for all $n\geq 2$ the scrolls $S(1^{n-2},2^2)$ and
$S(1^{n-1},3)$ can also be obtained by intersecting
$Seg(1,n+1)$ with a suitable $\Proj^{2n+1}$. \end{Lemma}

\begin{proof} By definition $X$ is isomorphic to $\p({\E})$,
with $\E$ a locally free sheaf of rank $n$ over a smooth curve
$C$. Since $X$ is regular, $C\simeq \Proj^1$ and hence
${\E}\simeq\oplus_{i=1}^n{\o}_{\Proj^1}(a_i)$, with $0<a_1\leq
a_2\leq\cdots a_n$. We know that $X$ is linearly normal and hence
that $2n+2=h^0(\Proj^1,{\E})=\sum_{i=1}^n(a_i+1)$ , from which it
follows
that $X$ is either $S(1^{n-2},2^2)$, $n\geq 2$,
or  $S(1^{n-1},3)$, $n\geq 1$. The
description as divisors in $Seg(1,n)$ as well as sections of
$Seg(1,n+1)$ is well
known. We leave the easy proof to the reader.\end{proof}

In view of Proposition \ref{linorm} and of Remark \ref {rationality}
(see also Proposition \ref {RPT} below), an {\it
OADP}-variety
$X\subset
\Proj^{2n+1}$ can be a scroll over a curve only if it is a smooth
rational normal scroll, i.e. only if $X$ is
either of type $S(1^{n-1},3)$ or of type $S(1^{n-2},2^2)$.\par

We will now prove, using an argument of F. Zak, that $S(1^{n-1},3)$ and
$S(1^{n-2},2^2)$ are indeed {\it OADP}-varieties for all $n\geq 2$.
Notice that $S(1^{n-1},3)$ makes sense even for
$n=1$, i.e. one has $S(3)$, a rational normal cubic, which is indeed an
{\it OADP}-variety.

\begin{Proposition}\label{OADPscrolls} For all $n\geq 2$ the scrolls
$S(1^{n-1},3)$ and
$S(1^{n-2},2^2)$ are {\it OADP}-varieties. These are the only scrolls
over a curve which are {\it
OADP}-varieties.\end{Proposition}

\begin {proof} By Lemma \ref {basiscrolls} the scrolls in question can
be obtained by intersecting $Y:=Seg(1,n+1)$
with a $\Proj^{2n+1}$. \par

Let now $p\in \Proj^{2n+3}$ be any point not on $Y$. It is well known
(see for instance \cite{Edge}) that the entry
locus $S_p(Y)$ has the form $\Proj^1\times\Proj^1_p$ for some
$\Proj^1_p\subset\Proj^{n+1}$, hence it is a
$2$-dimensional quadric in a linear $\Proj^3_p$. We are then in position
to apply Lemma \ref {zak} to
conclude with the first assertion.

The final assertion is a consequence of Lemma \ref {basiscrolls} and of
Proposition \ref{linorm} part iii) and Remark \ref{rationality}. \end{proof}
\end{ex}\bigskip

\begin{ex} \label{edgevarieties} ({\it Edge varieties}).
Consider $Seg(1,n)=S(1^{n+1})\subset \Proj^{2n+1}$ which has degree
$n+1$, $n\geq 1$. Recall that
$Pic(Seg(1,n))\simeq \Z<\Pi>\oplus \Z<\Sigma>$, where $\Pi$ is one
of the $\Proj^n$'s of the rulings of $S(1^{n+1})$ and $\Sigma$ is
$Seg(1,n-1)$ naturally embedded in $Seg(1,n)$.
Hence, an effective divisor $D$ of type $(a,b)$ has degree $a+nb$ in
$\Proj^{2n+1}$. The hyperplane divisor of
$Seg(1,n)$ is of type $(1,1)$.\par

Let $Q$ be a
general quadric of $\Proj^{2n+1}$ containing $\Pi$. Then the
residual intersection of $Q$ with $Seg(1,n)$ is an
$n$-dimensional smooth variety $E_{n,2n+1}$ of degree $2n+1$ in
$\Proj^{2n+1}$ which is a general divisor of type $(1,2)$ on $Seg(1,n)$.
Its degree is $2n+1$.

Let $\Pi_1$ and $\Pi_2$ be two fixed $\Proj^n$'s
of the ruling of $Seg(1,n)$. If $Q$ is a general quadric
through $\Pi_1$ and $\Pi_2$, then the residual
intersection of $Q$ with $Seg(1,n)$ is an $n$-dimensional smooth
variety $E_{n,2n}$ of degree $2n$ which is a general divisor of type
$(0,2)$ on $Seg(1,n)$. Notice that $E_{n,2n}\simeq\Proj^1\times Q'$,
where $Q'$ is a general quadric of $
\Proj^n$. In particular $E_{3,6}$ coincides with $Seg(1^3)$. \par

Varieties of the type $E_{n,2n+\epsilon}$, with $\epsilon=0,1$, have been
considered in \cite{Edge}. Edge proved
that they are {\it OADP}-varieties as soon as $n\geq 2$. Therefore they are
called {\it Edge varieties}. If $n=1$, then $E_{1,3}$ is the rational normal
cubic, which is an {\it OADP}-variety, whereas $E_{1,2}$ is a pair of
skew lines which, strictly
speaking, cannot be considered an {\it OADP}-variety because it is not
irreducible. However it still has one
apparent double point.\par

We remark that, by degree reasons, the Edge varieties and the scrolls
are different examples of {\it
OADP}-varieties, as soon as $n\geq 3$. If $n=2$ instead, $E_{2,4}$ is a
rational normal scroll of type $S(2^2)$. Notice that, instead, the other scroll
$S(1,3)$, which is also an {\it OADP}-surface, is {\it not} an  Edge
variety.\end{ex}\bigskip

We now sketch Edge's argument from \cite{Edge} to the effect that smooth
irreducible
divisors of type $(0,2)$, (1,2) and $(2,1)$ on the Segre
varieties $Seg(1,n)$, $n\geq 2$, have one apparent double
point. This generalizes the trivial fact that the only smooth
curves on a smooth quadric in $\Proj^3$ having one apparent
double point are of the above types.

\begin
{Proposition}\label{Le:Edge} Let $X\subset\Proj^{2n+1}$ be a smooth,
irreducible,
$n$-dimensional variety contained as a divisor of type $(a,b)$ in
$Seg(1,n)$,  $n\geq
2$. Then $X$ has one apparent double point if and only if
$(a,b)\in\{(2,1),(0,2), (1,2)\}$, i.e. if and only if $X$ is
either a rational normal scroll or one of the two Edge' s varieties as
in examples \ref {scrolls} and \ref
{edgevarieties}.  \end{Proposition}

\begin{proof} [Sketch of the proof.] As we know for, $p\notin
Y:=Seg(1,n)$  the  entry locus $S_p(Y)$ has the form
$Seg(1^2)\simeq \Proj^1\times\Proj^1_p$ for some
$\Proj^1_p\subset\Proj^n$ and spans a linear $\Proj^3_p$. So if
$X$ is a divisor of type $(a,b)$ of $Y$, the secant lines of $X$ passing
through $p$
are exactly the secant lines of $X\cap\Proj^3_p$ passing
through $p$. For a general $p\in\Proj^{2n+1}$, $X\cap\Proj^3_p$ is
a reduced, not necessarily irreducible curve and it is a divisor
of type $(a,b)$ on $\Proj^1\times\Proj^1_p$. Hence $X$ has one
apparent double point if and only if $(a,b)\in\{(1,2),\; (2,1),\;
(2,0),\; (0,2)\}$. If $(a,b)=(2,0)$, then
$X={\Proj}^n\amalg{\Proj}^n$ is reducible.
\end{proof}

In relation with Edge's varieties we notice the following useful
 generalization of \cite[Lemma 3 and Proposition 8]{Ru},
the main ingredients of the proof of Severi' s classification
Theorem of $OADP$-surfaces given in \cite{Ru}.

\begin{Proposition}
\label{Prop:hyperquadrics} Let $X\subset\p^{2n+1}$ be an {\it
OADP}-variety of dimension $n\geq 2$. Suppose $X$ contains a
$1$-dimensional family ${\cal D}$ of divisors $D$ of degree
$\alpha$, whose general element   spans  a $\p^n$. Then
$\alpha=2$, $X$ is an Edge variety, hence $X$ is either a divisor
of type $(0,2)$ or $(1,2)$ on $Seg(1,n)$, and ${\cal D}$ is cut
out on $X$ by divisors of type $(1,0)$.\end{Proposition}

\begin{proof} The assertion is true if $n=2$ (see \cite {Ru}, l.c.), so we
assume $n\geq 3$ from now on. By Zak's Theorem on Tangencies,
\cite[Corollary 1.8]{Zak1},  every  divisor $D\in\cal D$ has at
most a finite number of singular points, so that it is
irreducible.\par

First we will prove that ${\cal D}$ is a base point free pencil,
which is rational because $X$ is regular (see Remark \ref
{rationality}). In order to see this, it suffices to prove that
two general divisors $D_1,D_2$ of ${\cal D}$ do not intersect.
Suppose this is not the case. Then $D_1$ and $D_2$ intersect
along a variety $D_{1,2}$ of dimension $n-2$. Let $x$ be a
general point in a component of $D_{1,2}$ and let $p_i$ be a
general point of $D_i$, $i=1,2$. Of course $p_1,p_2$ are general
points of $X$. Notice that $\p^{n-2}\simeq T_xD_{1,2}\subset
<D_i>\simeq \p^n$, $i=1,2$. Since $T_{p_i}D_i\subset <D_i>$, then
$\dim(T_{p_i}D_i\cap T_xD_{1,2})\geq n-3$, $i=1,2$, hence $\dim
(T_{p_1}D_1\cap T_{p_2}D_2)\geq n-4$ and therefore also
$\dim(T_{p_1}X\cap T_{p_2}X)\geq n-4$. Since $X$ is not
defective, by Terracini's Lemma we have a contradiction as soon
as $n\geq 4$. If $n=3$, by arguing as above we see that
$<D_1>\cap <D_2>$ cannot be a plane.  Thus $D_{1,2}$ has to be a
line along which $D_1$ and $D_2$ intersect transversely.
 We
claim that $<x,p_i>\cap X=\{x,p_i\}$. Indeed, if $\alpha\geq 3$,
the line $<x,p_i>$ would intersect the divisor $D_i$ in a point
$p_i'$ different from $x$ and $p_i$ and the line $<p_1',p_2'>$
would cut the general secant line $<p_1,p_2>$ in a point
different from $p_1$ and $p_2$ and $X$ would not be an {\it
OADP}-variety. The claim proves that $\alpha=2$ and that $D$ is a
smooth quadric being smooth along a line contained in it. The
exact sequence $$0\to\o_X\to\o_X(D)\to\o_D(D)\to 0,$$ together
with $h^1(X,\o_X)=0$ and $h^0(D,\o_D(D))=2$, yields
$h^0(X,\o_X(D))=3$. The general hyperplane section $H$ of $X$
would be a smooth surface containing a $2$-dimensional family of
conics, and therefore $H$ would either be the Veronese surface in
$\Proj^5$ or a projection of it (e.g. see \cite{mez}), which is
not possible.\par

Now we consider the scroll $Y$ swept out by $<D>$ as $D$ varies
in ${\cal D}$. This is a  scroll over a rational curve
since there is a unique $\Proj^n$ through
the general point of $X$ and since $X$ is regular. Furthermore $X$ is  linearly normal, therefore
$Y$ is a rational normal scroll, i.e. a  cone over a smooth rational normal scroll
or $Seg(1,n)$ (see for instance
\cite[Theorem 2]{EH}).
 Let us prove that $Y$ cannot be a cone. Let
$V$ be the vertex of $Y$ and let $s\in V$ be any point. Keeping
the above notation, we remark that $s\not\in D_i$ and the lines
$<s,p_1>$ and $<s,p_2>$ span a plane $\Pi$, which contains the
general secant line $<p_1,p_2>$ to $X$. The line $<s,p_i>$
intersects the divisor $D_i$ at another point $p_i'$ different
from $p_i$, because it is a general line through $s\notin D_i$ in
the linear space $<D_i>$, $i=1,2$.  Then the secant line
$<p_1',p_2'>$ is contained in $\Pi$ and intersects the general
secant line $<p_1,p_2>$ outside $X$, implying that $X$ is not an
{\it OADP}-variety.\par

In conclusion, $Y$ is smooth, then $Y=Seg(1,n)$ and we conclude
by Proposition \ref {Le:Edge}.
\end{proof}

\begin{ex} \label {degree8} ({\it {\it OADP} $3$-folds of degree $8$ in
$\Proj^7$}). Let $X=\p({\E})\subset
\Proj^7$ be the scroll in lines over $\Proj^2$ associated to the very
ample vector bundle
${\E}$ of rank 2, given as an extension by the following exact
sequence

\begin{equation}\label{sequence}
0\to\O_{\Proj^2}\to{\E}\to\I_{\{p_1,\ldots,p_8\},{\Proj^2}}(4)\to 0, \end{equation}

\noindent where $p_1,\ldots,p_8$ are distinct points in $\Proj^2$, no
$7$ of them lying on
a conic and no $4$ of them collinear. Since $h^0(\Proj^2,{\E}(-2))=0$
such a vector bundle is stable, it has Chern classes
$c_1({\E})=4$, $c_2({\E})=8$ and it exists (see \cite{LP} and
\cite{Io2}).
 A general hyperplane section of $X$ corresponds to a section
 as in sequence (\ref{sequence}) and it is a smooth octic
rational surface $Y$, the embedding of the blow-up $\tilde \Proj^2$ of $\Proj^2$ at
the points $p_i$ given by the quartics through them.  The surface $Y$ is
arithmetically
Cohen-Macaulay and it is cut out by $7$ quadrics which define a Cremona
transformation $\varphi:\Proj^6\map\Proj^6$
(see \cite{ST2}, \cite{HKS}, \cite {MR}  and also \cite{AR}, \S 6).
Hence $X$ is an arithmetically Cohen-Macaulay 3-fold cut out by 7
quadrics and it is an {\it OADP}-variety by
Proposition \ref {cremona}. \par

We sketch now another description of $X$, from which we will deduce another proof
of the fact that $X$ is an {\it OADP}-variety. Notice that
$h^0(\Proj^2,I_{\{p_1,\ldots,p_8\},{\Proj^2}}(3))=2$ and call
$f_i:=f_i(x_0,x_1,x_2)$ two independent cubics through $p_1,\ldots,p_8$. Then
$H^0(\Proj^2,I_{\{p_1,\ldots,p_8\},{\Proj^2}}(4))$ is spanned by $f_ix_j, i=0,1,
j=0,1,2$, and an irreducible quartic, $\phi:=\phi(x_0,x_1,x_2)$,   vanishing at
$p_1,\ldots,p_8$. Then $H^0(\Proj^2,{\E})$ is spanned by sections
$s_{i,j}$ and $s$ respectively mapping to $f_ix_j$, $i=0,1, j=0,1,2$, and to
$\phi$ via the surjective map  $H^0(\Proj^2,{\E})\to
H^0(\Proj^2,I_{\{p_1,\ldots,p_8\},{\Proj^2}}(4))$ deduced by sequence (\ref{sequence}),  and by
the section $\sigma$ given by the inclusion $\O_{\Proj^2}\to{\E}$ also deduced from
(\ref{sequence}). Introduce homogeneous coordinates $z_{i,j}$, $i=0,1, j=0,1,2$,
$z_0,z_1$ in $\Proj^7$. Then the inclusion $X=\Proj({\E})\to \Proj^7$ is defined by

$$z_0=\sigma, z_1=s, z_{i,j}=s_{i,j}, i=0,1, j=0,1,2.$$

\noindent Whence we see that $X$ verifies the equations
$$z_{i,j}z_{h,k}=z_{i,k}z_{h,j}, i,h=0,1, j,k=0,1,2.$$
which, in the $\Proj^5$ with coordinates $z_{i,j}$, $i=0,1, j=0,1,2$,
define $Seg(1,2)$. Thus the line $L$ with equations $z_{i,j}=0, i=0,1, j=0,1,2$,
sits on $X$ and by projecting $X$ from $L$ to $\Proj^5$ we have $Seg(1,2)$. In
other words $X$ sits on the cone $F$ with vertex $L$ over $Seg(1,2)$. The variety
$X$ is the complete intersection of two divisors of type $(1,2)$ on $F$ (see \cite
{HKS}). Notice that a divisor of type $(1,2)$ on
$F$ is the intersection of $F$ with a quadric containing a $\Proj^4$ of
the ruling.\par

This description allows us to
see in another way that $X$ is an {\it OADP}-variety.
Indeed, let $p$ be a general point of
$\Proj^7$. Let $p'$ be its projection to $\Proj^5$ from $L$. The entry
locus $S_{p'}(Seg(1,2))$ is a
$2$-dimensional smooth quadric $Q$. Consider the $\Proj^5$ spanned by
$Q$ and $L$. Any secant to $X$ through $p$
sits in it.  The intersection of $X$ with this $\Proj^5$ is the
intersection of the rank $4$ quadric cone with
vertex $L$ over $Q$ with two divisors of type $(1,2)$. This is a
rational normal scroll of degree $4$, which is an
{\it OADP}-surface. This yields that there is a unique secant line to $X$
through $p$. \par

The focal locus of $X$ contains the cone $F$.\par\end{ex}\bigskip

The following examples have been first pointed out  by F. Zak.

\begin{ex} \label{spinor} {\it (A smooth $4$-fold of degree $12$).} For
every $k\geq 1$, we denote by
$S^{(k)}\subset \Proj^{2^k-1}$ the {\it $k$-th spinor variety} of
dimension $k(k+1)/2$,
parameterizing the subspaces of dimension $k-1$ contained in a smooth
quadric $Q_{2k- 1}\subset\Proj^{2k}$.\par

Consider $Y:=S^{(4)}\subset \Proj^{15}$. It is known that for
$p\in\Proj^{15}$ general, the entry locus
$S_p(Y)$ is a smooth $6$-dimensional quadric (see \cite {ES}). By lemma
\ref{zak}, the intersection $X$ of $Y$
with a general $\Proj^9$ is an {\it OADP}-variety of dimension $4$. \par

A different, but related, proof that $X$ is an {\it OADP}-variety, can
be obtained by applying Proposition
\ref {cremona} and by using the fact that the linear system $|{\cal
I}_{Y,\Proj^{15}}(2)|$ maps $\Proj^{15}$ to a
smooth quadric in $\Proj^9$, with fibers the $\Proj^7$'s spanned by the
entry loci $S_p(Y)$ (see \cite {AR}). \par

Notice that $Y$ is a {\it Mukai variety}, i.e. its general curve section
is a canonical curve (see \cite {Mu}). The
same happens for $X$. \end{ex}

\begin{ex} \label{lagrassmann} {\it (Lagrangian Grassmannians).} Let us
consider the {\it Lagrangian
Grassmannians} $\G^{lag}_k(2,5)$, over the four composition algebras
$k=\R,\C,\HH,\OO$, embedded via their
Pl\"ucker embedding.

These varieties can be described also  as
{\it twisted cubics
over  the four  cubic Jordan algebras} $Sym^3{\mathbb C}$, $M_3{\mathbb
C}$, $Alt_6{\mathbb C}$, respectively $H_3{\mathbb O}$, see \cite{Mu2}.
Note that also
Edge varieties $E_{n,2n}$ are of this type for split Jordan
algebras, see \cite{Mu2}.
Another possible description is as follows, see \cite{Zak2}:
\begin{itemize}
\item[i)] $\G^{lag}_\R(2,5)$ is the {\it symplectic Grassmannians}
$\G^{symp}(2,5)$ of the
$2$-planes in $\Proj^5$ on which a non-degenerate null-polarity is
trivial, i.e. which are isotropic with respect
to a non-degenerate antisymmetric form. This is a variety of dimension
$6$ which turns out to be embedded via the
Pl\"ucker embedding in $\Proj^{13}$,
\item[ii)]$\G^{lag}_\C(2,5)$ is the usual Grassmannians $\G(2,5)$
of $2$-planes in $\Proj^5$ which has
dimension $9$ and sits in $\Proj^{19}$ via the Pl\"ucker
embedding,
\item[iii)]$\G^{lag}_\HH (2,5)$ is the spinor variety $S^{(5)}$
of dimension $15$ in $\Proj^{31}$,
\item[iv)] $\G^{lag}_\OO(2,5)$ is the $E_7$-variety of dimension
$27$ in $\Proj^{55}$.
\end{itemize}

These are {\it OADP}-varieties. For the proof, we refer the reader to
\cite {AR}.  We notice that
$\G^{lag}_\R(2,5)$ is a Mukai variety (see \cite {Mu}). \end{ex}

\section{Degenerations of projections}\label {degenerations}

One of our tools for studying {\it OADP}- varieties, or, more generally,
{\it
RPT$_i$}-varieties, is the analysis of certain degenerations of
projections, which, we believe, can be useful
also in other contexts. This is what we present in this section.
\par

Let $X\subset \Proj^r$ be a smooth, irreducible, non-degenerate,
non-defective, projective variety of dimension $n$. Let $\Pi\subset
\Proj^r$ be a
general projective subspace of dimension $s$ with $r-2n-1\leq s\leq
r-n-2$.
Consider the projection morphism $X\to \Proj^ {r-s-1}$ of $X$ from
$\Pi$.
Let $x$ be a general point on $X$ and let
$T_xX$ be the projective tangent space to $X$ at $x$. Roughly speaking,
in this
section we will study how the projection $X\to \Proj^ {r-s-1}$ {\it
degenerates}
when its center
$\Pi$ tends to a general $s$-dimensional subspace
$\Pi_0$ containing $x$, i.e. such that $\Pi_0\cap X=\Pi_0 \cap
T_xX=\{x\}$. To be
more precise we want to describe the {\it limit} of the {\it double
point scheme} of
$X\to \Proj^ {r-s-1}$ in such a degeneration.
\par

Let us describe in more details the set up in which we work. We let
$T$ be a
general
$\Proj^{n+s+1}$ which is tangent to $X$ at $x$, i.e. $T$ is a general
$\Proj^{n+s+1}$ containing $T_xX$. Then we choose a general line $L$
inside $T$
containing $x$, and we also choose
$\Sigma$ a general $\Proj^{s-1}$ inside $T$. For every $t\in L$, we let
$\Pi_t$
be the span of $t$ and $\Sigma$. For $t\in L$ a general point, $\Pi_t$
is a general
$\Proj ^s$ in
$\Proj^r$. For a general
$t\in L$, we denote by $\pi_t: X\to \Proj^ {r-s-1}$ the projection
morphism of $X$
from $\Pi_t$. We want to study the limit of $\pi_t$ when $t$ tends to
$x$. \par

In order to do this, consider a neighborhood $U$ of $x$ in $L$ such that $\pi_t$ is
a morphism for all $t\in U-\{x\}$. We will fix coordinates on $L$ so that $x$ has
the coordinate $0$, thus we may identify $U$ with a disk
around $x=0$ in ${\C}$. Consider the products ${\cal X}=X\times
U$,
$\Proj_U= \Proj^ {r-s-1}\times U$. The projections $\pi_t$, for $t\in U
$, fit
together to give a rational map $\pi: {\cal X}\map \Proj_U$, which is
defined
everywhere except in the pair $(x,x)=(x,0)$. In order to extend it, we
have to blow up ${\cal X}$ at
$(x,0)$. Let $p: \tilde{\cal X}\to {\cal X}$ be this blow up and let
$Z\simeq \Proj
^n$ be the exceptional divisor. Looking at the obvious morphism $\phi:
\tilde{\cal
X}\to U$, we see that this is a flat family of varieties over $U$. The
fiber $X_t$
over a point $t\in U\setminus\{0\}$ is isomorphic to
$X$, whereas the fiber $X_0$ over $t=0$ is of the form $X_0=\tilde X\cup
Z$, where
$\tilde X\to X$ is the blow up of $X$ at $x$, and $\tilde X\cap Z=E$ is
the
exceptional divisor of this blow up, the intersection being transverse.
On
$\tilde{\cal X}$ the projections $\pi_t$, for $t\in U $, fit together to
give a
{\it morphism} $\tilde \pi: \tilde{{\cal X}}\to \Proj_U$. If we set ${\cal
H}=p^*(\pi^*({\cal O}_{\Proj_U}(1))$, then the map $\tilde \pi$ is
determined by the
line bundle ${\cal H}\otimes {\cal O}_{\Proj_U}(-Z)$. The restriction of
$\tilde \pi$ to a fiber $X_t$, with
$t\in  U\setminus\{0\}$, is the projection $\pi_t$. By abusing notation, we will
denote by $\pi_0$ the restriction
of
$\tilde \pi$ to $X_0$.
\par

We want to understand how $\pi_0$ acts on $X_0$. The restriction of
$\pi_0$ to $\tilde X$ is the projection of $X$ from the subspace
$\Pi_0$. Notice
that, since
$\Pi_0\cap X=\Pi_0 \cap T_xX=\{x\}$, this projection is not defined on
$X$ but it is
well defined on
$\tilde X$. \par

The main point of our analysis is to describe the action of $\pi_0$ on
$Z$. This is
the content of the following lemma:

\begin{Lemma} In the above setting, $\pi_0$ maps isomorphically $Z$ to
the
$n$-dimensional linear space $S$ which is the projection of
$T$ from
$\Pi_0$. \label{proiezione}
\end{Lemma}

\begin{proof} The restriction of $\pi_0$ to $Z\simeq \Proj ^n$ is
determined by the
line bundle
${\cal O}_Z(-Z)\simeq {\cal O}_Z(E)\simeq {\cal O}_{\Proj ^n}(1)$. Since
$\pi_0$ is a
morphism, the image of $Z$ via $\pi_0$ is a $\Proj ^n$.\par

The projection $S$ of $T$ in $\Proj^ {r-s-1}$ from $\Pi_0$ has also
dimension $n$. Notice that $S$ contains the image $\Omega$, under the
same
projection, of $T_xX$. This is a subspace of dimension $n-1$ which is
also the image
of the exceptional divisor $E$ of
$\tilde X$ via $\pi_0$, and therefore it is contained in $\pi_0(Z)$ too.
\par

Consider the smooth curve $C$ in ${\cal X}$, whose general point is
$(x,t)$, with
$t\in U$ general, and let
$\tilde C$ be its proper transform on $\tilde {\cal X}$. We let $\gamma$
be the
point $\tilde C\cap Z$. \par

This curve is contracted by the obvious projection map $\tilde {\cal
X}\to
\Proj ^{r-s-1}$, to a point $p\in  \Proj ^{r-s-1}$. Of course
$p=\pi_0(\gamma)\in
\pi_0(Z)$. Furthermore $p\in S$. This follows from the fact that $p\in
S_t$, for all
$t\in U\setminus\{x\}$, where $S_t$ is the projection of $T$ via $\pi_t$. We
claim that $p\notin \Omega$. Otherwise the $\Proj^{s+1}$ joining $L$ and
$\Sigma$
would cut $T_xX$ along a line, against the generality assumptions about
$L$ and
$\Sigma$. Finally our assertion follows: indeed $\pi_0(Z)$ and $S$ are
two
$\Proj^n$'s, meeting along $\Omega$ which is a $\Proj^{n-1}$, and also
meeting at
a point
$p\notin \Omega$. \end{proof}

The {\it double point scheme} of the map $\tilde \pi$ is the scheme
$\tilde{\cal X}\times _{\Proj_U}\tilde{\cal X}$. Consider the obvious
map $\psi:
\tilde{\cal X}\times _{\Proj_U}\tilde{\cal X}\to U$. We will denote by
$\Delta_t$
the fiber of $\psi$ over $t\in U$. We will assume, up to shrinking $U$,
that $\dim(\Delta_t)$ is the minimum
possible, i.e.

$$\dim (\Delta_t)=2n-r+s+1\geq 0$$

\noindent for all $t\in U\setminus\{0\}$. The family $\psi: \tilde{\cal X}
\times_{\Proj_U}\tilde{\cal X}\to U$ may very well be not flat over
$x=0$. However, there
is a unique flat limit of
$\Delta_t$,
$t\not= 0$, sitting inside $\Delta_0$. We will denote by $\tilde
\Delta_0$ this flat limit, and we will call it
the {\it limit double point scheme} of the projection $\pi_t$, $t\not=
0$. We want to
give some information about it. We notice the following lemma:

\begin{Lemma}  In the above setting, every irreducible component of
$\Delta_0$ of
dimension
$2n-r+s+1$ sits in the limit double point scheme $\tilde \Delta_0$.
\label{flatlimit}
\end{Lemma}

\begin{proof} Every irreducible component of
$\tilde{\cal X}\times _{\Proj_U}\tilde{\cal X}$ has dimension at least
$2n-r+s+2$. Hence every
irreducible component of
$\Delta_0$ of dimension $2n-r+s+1$ sits in a component of  $\tilde{\cal
X}\times _{\Proj_U}\tilde{\cal X}$ of dimension exactly $2n-r+s+2$
dominating $U$. The assertion immediately
follows. \end{proof}

Let us introduce some more notation. We will denote by:
\begin{itemize}
\item[i)] $X_T$ the scheme cut out by
$T$ on $X$. $X_T$ is cut out on $X$ by $r-n-s-1$ general
hyperplanes
tangents at
$x$. We call $X_T$ a {\it general \hbox{$(r-n-s-1)$-tangent} section}.
Remark that
each component of
$X_T$ has dimension at least $2n-r+s+1$,
\item[ii)] $\tilde X_T$ the proper transform of $X_T$ on $\tilde
X$,
\item[iii)] $Y_T$ the image of $\tilde X_T$ via the restriction of
$\pi_0$ to
$\tilde X$. By Lemma \ref{proiezione}, $Y_T$ sits in
$S=\pi_0(Z)$, which is naturally isomorphic to $Z$. Hence we may
consider $Y_T$
as a subscheme of $Z$,
\item[iv)]$\Delta'_0$ the double point scheme of the restriction of
$\pi_0$ to
$\tilde X$.
\end{itemize}

With this notation in mind, the next lemma follows right away:

\begin{Lemma} In the above setting, $\Delta_0$ is the union of
$\Delta_0'$ and
$\tilde X_T$ on $\tilde X$ and of $Y_T$ on $Z$.
\label{doublepointscheme}
\end{Lemma}

As an immediate consequence of Lemma \ref{flatlimit} and Lemma \ref{doublepointscheme},
we have the following proposition, which will play a
crucial role in what follows:

\begin{Proposition} In the above setting, every irreducible component of $X_T$ of
dimension $2n-r+s+1$ gives rise to an irreducible component of
$\tilde X_T$ which is contained in the limit double point scheme $\tilde
\Delta_0$. In particular if $X_T$ has
dimension $2n-r+s+1$, then
$\tilde X_T$ is contained in $\tilde\Delta_0$.
\label{limitdoublepointscheme}
\end{Proposition}\bigskip

\section{Applications to {\it RPT} and {\it OADP}-varieties} \label {applications}

One may apply the results from the previous section to the study of
$n$-dimensional
{\it RPT} and {\it OADP}-varieties. We start with the following:

\begin{Theorem}\label {prebronowski} Let $X\subset \Proj^{2n+1}$ be a smooth,
irreducible, projective variety of dimension $n$ whose number of
apparent double points is $\nu(X)>0$. Let $\delta(X)$ be the degree of the
projection of $X$ to $\Proj^n$ from the tangent $\Proj^n$ to $X$ at its general
point. Then $0<\delta(X)\leq \nu(X)$.\end {Theorem}

\begin{proof} Notice that $X$ is not defective, because $\nu(X)>0$. Hence
by Terracini's lemma $\delta(X)>0$ (see \cite {cc} or \cite {AR}, Proposition 3).
Now we apply Proposition \ref {limitdoublepointscheme} to $X$ with $s=0$. Notice
that, by assumption, $X_T$ has $\delta(X)$ isolated points, which give rise to
$\delta(X)$ points in the limit double point scheme. Since this has degree
$\nu(X)$, we have the assertion.\end{proof}

As a consequence we have the following.

\begin{Corollary} \label{bronowskiconj} Let $X\subset \Proj^{2n+1}$ be a smooth,
irreducible, projective {\it OADP}-variety of dimension $n$. Then the projection of
$X$ to $\Proj^n$ from the tangent $\Proj^n$ to $X$ at its general point is
birational. In particular $X$ is rational.\end{Corollary}

\begin{remark} \label {Bronowski}  J. Bronowski
claims in \cite {Br} that  $X_n\subset \Proj^{2n+1}$ is an {\it
OADP}-variety if and only if the projection of $X$ to $\Proj^n$ from a
general tangent $\Proj^n$ to $X$ is a
birational map. Unfortunately Bronowski's argument is very obscure and,
at the best of our knowledge, there is no convincing proof in the present
literature for this statement as it is. It is therefore sensible to consider it as
a conjecture. The above corollary proves one implication of Bronowski's
conjecture. \par

For more information about this  conjecture see \cite {AR}.
\end{remark}
\medskip

Before stating and proving the next result, we need to recall a basic fact concerning degenerations of
curves (see also \cite {hi}):

\begin{Proposition} \label{classical} Let $p: {\cal C}\to \Delta$ be a
proper, flat family of curves parametrized
by a disk $\Delta$. Suppose that the general curve $C_t, t\in
\Delta\setminus\{0\}$, of the family is irreducible of
geometric genus $g$. Then every irreducible component of the central
fiber $C_0$ has geometric genus $g'\leq
g$.\end {Proposition}

\begin {proof} Let $f: {\cal C}'\to {\cal C}$ be the normalization map
and set $p'=p\circ f$. Then
$p': {\cal C}'\to \Delta$ is again a flat family and its general curve
$C'_t, t\in \Delta\setminus\{0\}$, is smooth of
genus $g$. Moreover the central fiber $C'_0$ has no embedded points and
it is a partial normalization of $C_0$.
\par

Let $p'': {\cal C}''\to \Delta$ be a semistable reduction of $p': {\cal
C}'\to \Delta$ such that all the
irreducible components of the central fiber $C''_0$ are smooth. The
curve $C''_0$ is reduced, connected (see \cite[Exercise 11.4, pg. 281]{ha}),
of arithmetic genus $g$. Hence all its
irreducible components have geometric genus
bounded above by $g$. Every irreducible component of $C_0$ is birational
to some
irreducible component of $C'_0$, which, in turn, is dominated by some
irreducible component of $C''_0$. This
proves the assertion.  \end{proof}

\begin {Corollary}\label{uniruled} Let $p: {\cal F}\to \Delta$ be a
proper, flat family of varieties over a disk
$\Delta$ such that the general fiber is smooth, irreducible and
uniruled. Then every irreducible component of the
central fiber is also uniruled.\end{Corollary}

\begin{proof} There is no lack of generality in assuming that $\cal F$ is
smooth. Take an
irreducible component of the central fiber and let $x$ be a general
point on it. We can consider a section of
the family, i.e. a morphism $\gamma: U\to {\cal F}$ such that
$p(\gamma(t))=t$, for all $t\in U$, and such that
$\gamma(0)=x$. As a consequence of the uniruledness assumption, there
is, up to a harmless base change, a
subvariety ${\cal C}$ of ${\cal F}$ such that the restriction $p: {\cal
C}\to \Delta\setminus\{0\}$ is a flat family of
rational curves such that for all $t\in  \Delta\setminus\{0\}$ the curve $C_t$
of the family contains $\gamma(t)$. The
flat limit $C_0$ of this family of curves contains $x$ and all of its
components are rational by Proposition \ref{classical}. This proves the assertion.\end {proof}

Now we are ready for the announced result about {\it RPT}-varieties.

\begin{Proposition} \label{RPT} Let $X\subset \Proj^r$ be a smooth,
irreducible, non-degenerate, projective variety
of dimension $n$.
\begin{itemize}
\item[i)] Suppose $X$ is a {\it RPT}-variety and let $C$ be its
general
$(n-1)$-tangent section. Then every
$1$-dimensional irreducible component of $C$ is rational.
\item[ii)]Suppose $X$ is a {\it RPT$_{n-1}$}-variety of dimension
$n\geq 3$,
and let $D$ be its general tangent hyperplane section. Assume $X$ is not
a scroll over a curve. Then $D$ is
irreducible and uniruled and the general hyperplane section
$H$ of $X$ has $\kappa(H)=-\infty$. Furthermore $X$ is also uniruled,
hence $\kappa(X)=-\infty$.
Finally $h^i(X,{\cal O}_X)=0$, $i=n-1,n$.
\end{itemize}
\end{Proposition}

\begin{proof} The first assertion of $(i)$ is an immediate consequence
of
Proposition \ref {limitdoublepointscheme} and Proposition \ref {classical}. \par

As for $(ii)$, if $X$ is not a scroll, its general tangent hyperplane
section $D$
has double points along a subspace $\Pi_D$ of dimension $s$ such
that $0\leq
s\leq n-3$ and it is smooth
elsewhere (see \cite {Kl} and \cite[Theorems 2.1 and 3.2]{Ei}). This implies that $D$,
which is
connected, is also irreducible (see also Proposition \ref{OADPrigate}
below). By  Proposition
\ref {limitdoublepointscheme} and Corollary \ref {uniruled},
$D$ is uniruled. As a consequence $X$ is uniruled and
$\kappa(X)=-\infty$.\par

Consider a
flat family ${\cal F}$ over a disk $U$ having $H$ as the general fiber
and $D$ as the
central fiber. By semicontinuity we have
$h^0(D,K_{{\cal F}|D}^{\otimes i})\geq h^0(H,K_{{\cal F}|H}^{\otimes
i})=h^0(H,K_H^{\otimes i})$ for all $i\geq 1$. Let $\tilde D\to D$ be a
desingularization of
$D$. Since
$D$ is uniruled and has canonical singularities, then $h^0(D,K_{{\cal
F}|D}^{\otimes
i})= h^0( D,K_ D^{\otimes i})=h^0(\tilde D,K_{\tilde D}^{\otimes
i})=0$ for all $i\geq 1$. This implies
$\kappa(H)=-\infty$. \par

>From the exact sequence

$$0\to {\cal O}_X(-1)\to {\cal O}_X\to {\cal O}_H\to 0$$

\noindent and by Kodaira vanishing theorem, we see that $h^{n-1}(X,{\cal
O}_X)\leq h^{n-1}(H, {\cal
O}_H)=h^0(H,{\cal O}_H(K_H))=0$. Moreover $h^n(X,{\cal O}_X)= h^0(X,
{\cal O}_X(K_X))=0$\end{proof}

In the surface case we can prove a slightly stronger result:

\begin{Proposition} \label {PGsurfaces} Let $X\subset \Proj^r$, $r\geq
4$, be a
smooth, irreducible, non-degenerate surface which is different from the
Veronese
surface in $\Proj^r$, $r=4,5$. Then $PG(X)$ is defined. Furthermore the
sectional
genus $\pi$ of $S$ satisfies $\pi \leq PG(X)+\mu$, where $\mu=0$ if $X$ is a
scroll, and
$\mu=1$ otherwise. \end{Proposition}

\begin{proof} The Veronese surface in $\Proj^5$ is the only smooth
defective
surface and the same surface, and its general projection in $\Proj^4$,
are the only
surfaces such that the double curve of their general projection in
$\Proj^3$ is
reducible (see
\cite{Se}, \cite{fr1}, \cite{fr2}).\par

The general tangent curve section $C$ of $X$ has a node at the contact
point and it
is smooth elsewhere. If $X$ is not a scroll, $C$ is also irreducible (cfr. Proposition \ref{scroll}). If
$X$ is a
scroll it consists of a general line $R$ of the ruling plus a unisecant
$C'$ meeting
$L$ at the contact point.  \par

By Proposition \ref {limitdoublepointscheme} and Proposition
\ref{classical} every irreducible component of $C$ has
genus $\gamma\leq PG(X)$. Hence, if $X$ is not a scroll, $C$
is irreducible of geometric genus
$\gamma\leq PG(X)$. Since $C$ has a node, its arithmetic genus is
$\gamma+1\leq
PG(X)+1$. On the other hand, the arithmetic genus of $C$ is the
sectional genus $\pi$ of $X$. The argument is
similar if $X$ is a scroll. \end{proof}

As a consequence, by taking into account the results from \cite {ab}, \cite{Io1}
and \cite {Io4} (see also \cite {li}), one may classify all surfaces $X$ with
$PG(X)\leq 4$.
We will not dwell on this now, and we will simply classify here the
{\it RPT}-surfaces:

\begin{Proposition} \label{RPTsurfaces} Let $X\subset \Proj^r$, $r\geq
4$, be a
smooth, irreducible, non-degenerate surface which is different from the
Veronese
surface in $\Proj^r$, $r=4,5$. Then $X$ is a {\it RPT}-surface if and
only if it is
one of the following:
\begin{itemize}
\item[i)] a rational normal scroll of degree 3 in $\Proj^4$,
\item[ii)] a rational normal scroll of degree $4$ in
$\Proj^5$,
\item[iii)] a del Pezzo surface of degree $n=4,5$ in $\Proj^n$.
\end{itemize}
\end{Proposition}

\begin{proof} By Proposition \ref {PGsurfaces}, the sectional genus of
$X$ is
$\pi\leq 1$ and $X$ cannot be an elliptic scroll. Hence $X$ is either a rational scroll or a del
Pezzo
surface. By checking the various possible cases and using the classical
formula for
the geometric genus of the double curve of the general projection of a
surface (see \cite[pg. 176]{enr} and also \cite {Pi}), one proves the
assertion.\end{proof}

\begin{remark} \label{p4} The hypothesis of smoothness for $X$ in the
previous
proposition can be relaxed. It suffices to assume that $X$ has at worse
improper double points.
\end {remark}
\medskip

By Proposition \ref {OADP-FRPT} all {\it OADP}-surfaces are {\it
RPT}-surfaces, i.e
they are listed in the previous proposition. Actually all surfaces there
lying in
$\Proj^5$ are  {\it OADP}-surfaces, and one finds a new, rather cheap,
proof
of the classification of {\it OADP}-surfaces (see \cite{Se}, \cite{Ru}):

\begin{Theorem} \label {OAPDsurfaces} Let $X\subset \Proj^5$ be an {\it
OADP}-surface. Then $X$ is one of the following:\medskip

\noindent (i) a rational normal scroll of degree $4$;\medskip

\noindent (ii) a del Pezzo surface of degree 5.\end {Theorem}\bigskip

\begin{remark}\label{rationalhypsect} We notice that an equally easy proof of
Theorem \ref {OAPDsurfaces} could be deduced from Corollary \ref
{bronowskiconj}. In particular for smooth non degenerate surfaces in $\Proj^5$,
OADP is equivalent to FRPT. It is reasonable to suspect that this relation holds also in
higher dimensions, see Remark \ref{re:OADPvsFRPT}.
\end{remark}

Turning to the classification of
{\it OADP}-threefolds $X$, which is our main aim, we notice that Proposition \ref {RPT} and Remark \ref
{rationality} imply that the general hyperplane  section of
such an $X$ is rational. A proof of the classification Theorem \ref
{Th:classification} below could be based
on this remark. Indeed by applying Mori's program (see \cite {mori}) as in \cite
{Io3}, \cite {Io5} or \cite {campflen}, one sees that $X$ is of one of the
following types (see \cite[Theorem II]{Io3} ):\begin{itemize}
\item[i)] a scroll over a curve or over a surface,
\item[ii)] a quadric fibration,
\item[iii)] a del Pezzo threefold,
\item[iv)] a quadric in $\Proj^4$,
\item[v)] (a projection of) a Veronese surface fibration,
\item[vi)] a projection of the Veronese embedding of $\Proj^3$ via the cubics
or of the Veronese embedding of a quadric in $\Proj^4$ via the quadrics.
\end{itemize}

A case by case analysis leads to the classification. Most of the tools
needed for this have already been
introduced so far or will be in \S \ref {tanbeh}. However, the
approach we have chosen to follow here
 is slightly different and, we think, somewhat conceptually easier. It
is,
in
a sense, dual to the one outlined above. In fact one first bounds the
degree of $X$, drastically reducing the a priori
possible cases (see next \S \ref{tanbeh}: this is not more expensive than
the analysis required in the discussion of
the cases $(i),...,(vi)$ above). Then a bit of Mori's theory, i.e.
adjunction theory, is used to determine, for each
possible degree, the varieties that really occur.

\section{Tangential behavior of {\it OADP} 3-folds}\label {tanbeh}

Let $X\subset\Proj^7$ be an {\it OADP}-variety of dimension $3$. The
main result of this
section is an upper bound for the degree $d$ of
$X$.

\begin{Theorem} \label{gradoOADP} Let $X\subset\Proj^7$ be a {\it OADP}
$3$-fold of degree $d$. Then $d\leq 9$. Furthermore either  $X$ is a scroll over a curve and
$d=5$ or
$d=9-k$, where $k\leq 3$
is the number of lines passing through a general point of $X$
\end{Theorem}

The proof is based on the analysis of the {\it tangential
behavior} of $X$, i.e. on the study
of the intersection scheme $Z:=Z_x$ of $X$ with $T_xX$, where $x\in X$
is a general point. We  need a
few technical results, some of them of independent interest. We state them at the first occurrence,
but we will prove them separately  in
the next section, in order not to break the main line of the argument.
We refer to them freely here.\par

First we remark
that $0\leq \dim(Z_x)\leq 2$. The case of $\dim(Z_x)=2$ is classically known.
\begin{Proposition}\label{scroll} Let $X\subset \Proj^r$ be an
irreducible, projective, non degenerate
variety of dimension $n$ with $1<n<r$. Let $x$ be a
general point of $X$. The intersection scheme of $X$ with $T_xX$ has an
irreducible
component of dimension $n-1$ if and only if either $X$ is a hypersurface
in $\Proj^{n+1}$ or $X$ is a scroll over
a curve.\end{Proposition}
We will prove this Proposition in the next section for the lack of a modern reference.
This allows to restrict our attention to the case of an {\it OADP} $3$-fold
$X\subset\Proj^7$ such that the scheme $Z_x$, with
$x\in X$ a general point, has dimension at most $1$. Equivalently we
may and will assume from now on that $X$ is
not a scroll.
The main difficulty will be to treat the one dimensional case.
A great part of this section is indeed devoted to prove the following Lemma

\begin{Lemma} If $\dim(Z_x)=1$,
then $(Z_x)_{red}$ is a bunch of $k$ lines through $x$.
\end{Lemma}

Let us now fix some notations which we are going to use throughout this
section.\par

\noindent Let $f:Y\ra X$ be the blow up of $X$ at $x$ with exceptional
divisor $E$. We let $H'$ be the pull back on $Y$ via $f$ of the general
hyperplane section $H$ of $X$.  We
consider the proper transform $\Lambda$ on $Y$ of the linear system of
tangent hyperplanes to $X$ at $x$, i.e.
$\Lambda=|H'-2E|$. Let $S$ be a general element in $\Lambda$, i.e. $S$
is the proper transform on $Y$ of a
general tangent hyperplane section $\Sigma$ of $X$. The surface $\Sigma$ has an
ordinary double point at $x$ and no other singularity, see for instance
\cite[Theorem 2.1]{Ei}. Thus $S$ is a
smooth surface and
the exceptional curve $F=E_{|S}$ is an irreducible $(-2)$-curve. \par

We will consider the {\it characteristic system} $\cal C$ of $\Lambda$
on $S$, i.e. the linear system cut out by
$\Lambda$ on $S$. We may write ${\cal C}={\cal
Z}+aF+{\cal L}$ where ${\cal L}$ is the movable part of $\cal C$ and
${\cal Z}+aF$ is the fixed part, with ${\cal
Z}$ not containing $F$. Since $X$ has non-degenerate tangential variety,
by a classical result of Terracini's, one
has $a=0$ (see \cite{terracini}, or \cite[pg. 417]{GH}, or also
\cite[Theorem 7.3 i)]{Lan}). Moreover the image
of $\cal Z$ via $f$ is the support of the $1$-dimensional part of $Z$.
\par

By abusing notation, we often denote by ${\cal C}$ and ${\cal L}$
also a curve of the linear systems in
question. Notice then the basic linear equivalence relation
$$H^{\prime}_{|S}\equiv {\cal C}+2F={\cal Z}+2F+{\cal L}$$
which we
frequently use later on.\par

Next we are going to study the linear system $\cal C$, by first studying
its movable part $\cal L$ and then its
fixed part $\cal Z$.

\begin{Proposition} \label{Le:br} In the above setting, $\L$ is a
complete
linear system which is  base point free and not composed with a pencil.
Moreover
$\L^2=1$ and $h^1(S,{{\cal O}_S({\cal  L})})=0$. In particular either
$Z_x$ has dimension $0$ supported at $x$ or
it is equidimensional of dimension $1$ off $x$. Moreover the morphism
$\phi_{\cal L}: S\to \Proj^2$ is a sequence
of blowing-ups. \end{Proposition}

\begin{proof} The threefold $X$  is linearly normal and, as we saw in remark
\ref{mattuck}, $h^1(X,\O_X)=0$.
Thus any hyperplane
section of $X$ is linearly normal. This means that the linear system
${\cal H}'$ cut out by $|H'|$ on $S$ is
complete. Notice that ${\cal C}={\cal H}'(-2F)$, which proves that $\cal C$, and therefore also $\cal L$, is
complete. \par

Since $X$ is non-defective, the projection of $X$
to $\Proj^3$ from the tangent space
$T_xX$ is dominant (see \cite {cc} or \cite[Proposition 3]{AR}, ). This
implies that the rational map $\phi_{\cal
L}: S\to \Proj^2$ determined by the linear system $\L$ is dominant,
namely $\cal L$ is not composed with a pencil.

By Proposition \ref{RPT}, the general curve in $\L$, which is
irreducible, is also rational.
Here we need the following auxiliary result, whose proof is postponed to the next section.

\begin {Proposition} \label{rational} Let $S$ be a smooth, irreducible,
projective, regular surface and let $\cal
M$ be a line bundle on $S$. Suppose that the general curve $C\in |\cal
M|$ is irreducible with
geometric genus $g$ and that $h^0(S,{\cal M})\geq g+2$. Then the
complete linear
system $|\cal M|$ is base point free. Furthermore $S$ is rational.
\end {Proposition}

By Proposition \ref{rational}, $\L$ is base point free. Therefore the general element in
$\L$ is a smooth rational curve. Hence
$K_{S}\cdot\L=-2-\L^2$, and, by Riemann--Roch formula, $\chi({{\cal
O}_S({\cal  L})})=2+\L^2$. On the other hand
one has $h^0(S,{{\cal O}_S({\cal  L})})=3$ and $h^2(S,{{\cal O}_S({\cal
L})})=h^0(S,{{\cal O}_S(K_{S}-{\cal
L})})=0$, because $\L$ is effective and $S$ is rational. Since
${\L}^2>0$ we finally obtain $h^1(S,{{\cal
O}_S({\cal  L})})=0$ and $\L^2=1$.\par

The rest of the assertion is trivial.\end{proof}

\begin{remark}\label{re:OADPvsFRPT}
 Notice that, as a consequence of the above proposition, one finds a new
proof, in
the threefold case, of Corollary \ref {bronowskiconj}. This also proves the same statement
for linearly normal {\it FRPT} 3-folds embedded in $\Proj^7$. Recall that Bronowski, \cite{Br},
claimed that {\it OADP} varieties are characterized by this property.
Summing all one can conjecture that OADP and {\it FRPT} are equivalent also, keep in mind Remark
\ref{rationalhypsect}, for
smooth non degenerate 3-folds in $\Proj^7$ and maybe a similar statement is true in any dimension.
\end{remark}
\medskip

The next step toward the proof of Theorem \ref {gradoOADP} is to study
the fixed
part $\cal Z$ of the system $\cal C$, i.e. the scheme $Z:=Z_x$. We will denote by
$Z_s$ the Zariski closure of the support of the scheme $Z_x\setminus\{x\}$. By Proposition
\ref {Le:br}, either $Z_s=\emptyset $ or every irreducible component of
$Z_s$ has dimension $1$. We will consider the latter case next, for this we use
the notion of $k$-filling curve.

\begin{Definition} Let $C\subset X$ be a reduced irreducible
curve on $X$. We say that $C$ is a {\it $k$-filling curve} for $X$ if $C$
is the general member of a
$k$-dimensional family of curves such that there is a curve of the
family passing through the general point of
$X$.
\end{Definition}

\begin{Lemma} If $\dim(Z)=1$, let $D\subset Z_s$ be any
irreducible component. Then $D$ is a $k$-filling curve for $X$ for some $k\geq 2$.
In particular if $Z_s$ contains a line then $x$ sits on some line contained in
$Z_s$. \label{le:kfil} \end{Lemma}

\begin{proof} Consider the family $\cal D$ in which $D$ varies as $x$ varies on
$X$. By abusing notation, we denote by $\cal D$ also the parameter space of
the family. \par

First we claim that $\dim({\cal D}):=k\geq 2$. Otherwise there
would be a  surface $\Delta$ on $X$ such that $D\subset T_xX$ for the general
point $x\in \Delta$. Then $1\leq \dim(<D>)\leq 2$. Suppose $D$ is
a line.  By generically projecting to $\Proj^4$ we would find a
hypersurface $X'$, the projection of $X$, containing
a line $D'$, the projection of $D$, such that the general hyperplane
through $D'$ is tangent to $X'$ at some
smooth point $x'\notin D'$. This contradicts Bertini's theorem. Also if $D$ spans
a plane a similar argument, which we can leave to the reader, works to reach a
contradiction to Bertini's theorem. This proves our assertion on $k$. \par

Now assume that all the curves $D$ of $\cal D$ are contained in a surface $B$.
Then through the general point of $B$ there is a $(k-1)$-dimensional family of
curves of $\cal D$. In particular two general curves  of $\cal D$ have to
intersect. This would imply that $T_xX\cap T_yX\neq\emptyset$ for $x$
and $y$ general in $X$, contradicting
Terracini's lemma. \end{proof}

\begin{Lemma} \label{Z1} If $\dim(Z)=1$, then either $Z_s$
consists of a bunch of lines passing through $x$ and spanning $T_xX$ or
it is a plane curve. In the latter case the plane spanned by $Z_s$ contains $x$.
\end{Lemma}

\begin{proof}  Recall that the focal locus $F(X)$ of $X$ has dimension
at most $5$ (see \S \ref {definitions}).
Thus $F(X)$ cannot contain the tangential variety $Tan(X)$, which has
dimension $6$. Therefore for the
general point of $z\in Tan(X)$, with $z\in T_xX$, $x\in X$, there is a
unique {\it secant} to $X$ through it, i.e.
the tangent line $<x,z>$. Since
$Z_s$ is a subvariety of $T_xX\simeq \Proj^3$, we
conclude that
$Sec(Z_s)$ is properly contained in $T_xX$.
In particular every irreducible component of $Z_s$ is planar. Moreover
the classical trisecant lemma (see, for
instance, \cite {cc}) yields that either $Z_s$ itself is planar or $Z_s$
consists of
a bunch of lines spanning $T_xX$ and passing through a point $y\in
T_xX$. In this latter case, since the general
tangent hyperplane section $\Sigma$ of $X$ contains $Z$, we conclude
that $x=y$, otherwise $\Sigma$ would be
singular at $y$, which, as we know, is impossible (see \cite[Theorems. 2.1]{Ei}).

If $Z_s$ is one single line, then $Z_s$ contains $x$ by Lemma \ref
{le:kfil}. Suppose finally $Z_s$ is a plane curve which is not a line and let
$\alpha$ be the plane where $Z_s$ sits. Thus $Sec(Z_s)=\alpha$. Suppose $x\notin
\alpha$. Then any tangent line $L$ at $x$ intersects $\alpha$ at a point
$y\not=x$ and all the secants to $Z_s$ through $y$ meet $L$. Thus $L$
would be focal, hence $F(X)$ would contain
$Tan(X)$, a contradiction.\end{proof}

Next assume that $Z_s$ is a plane curve. Let
$\alpha:=\alpha_x:=<Z_s>$ be a plane. We denote by $\delta\geq 2$ the
degree of the plane curve $Z_s$. We can consider the subvariety ${\cal
F}\in \G(2,7)$
described by the planes $\alpha_x$ when $x$ varies on $X$.

To better understand this variety we need the following Proposition,
whose proof is demanded.

\begin{Proposition}\label{tangentplanes} Let $X\subset \Proj^r$, $r\geq
4$, be a smooth, irreducible, non
degenerate, projective $3$-fold which is neither a scroll over a curve
nor a quadric in $\Proj^4$. Then there is
no irreducible subvariety $\cal F\subset \G(2,r)$ of dimension $2$,
such that the general point $p\in \cal F$
corresponds to a plane $\Pi\subset \Proj^r$ which is tangent to $X$
along an irreducible curve $\Gamma$
which is $2$-filling for $X$.\end{Proposition}

We are in the condition to prove the following lemmas.

\begin{Lemma}\label{tridimensional} In the above setting, one has
$\dim({\cal F})=3$. \end{Lemma}

\begin{proof}  Suppose the assertion is not true. Given a general plane
$\alpha\in {\cal F}$,
there would be an irreducible positive dimensional subvariety $\Delta$
of $X$ such that $\alpha$ would correspond
to the general point $x\in \Delta$. By Lemma \ref {Z1}, $\Delta$
has to lie on $\alpha$.  By Proposition
 \ref{scroll} and since we are assuming that $X$ is not a scroll, we have
that $\Delta$ is a curve. In conclusion the
planes $\alpha$ would vary in a $2$-dimensional family $\cal F$ and each
of them would be tangent to $X$ along a
curve $\Delta$, contrary to Proposition \ref{tangentplanes}.
\end{proof}

\begin{Lemma} \label {fuochi} In the above setting, one has $\delta\leq
3$.\end{Lemma}

\begin{proof}  Assume by contradiction that  $\delta\geq 4$.  By Lemma \ref {tridimensional},
and by applying the results of
\cite[\S 3]{mez} to a general projection of $X$ to
$\Proj^5$, we see that $X$ sits on a $4$-dimensional scroll $W$ over a curve.\par

Let $\Pi$ be the general $\Proj^3$ of the ruling of $W$ and let $D$ be the
surface cut out by $\Pi$ on $X$. Notice that $D$ is smooth in codimension $1$,
otherwise $\Pi$ would be tangent to $X$ along a curve, against Zak's theorem on
tangencies (see \cite {Zak1}). In particular $D$ is irreducible. Then by
Proposition \ref {Prop:hyperquadrics} we have that $X$ should be an Edge variety,
which has no plane curves of degree $\delta\geq 3$, a contradiction. \end{proof}

\begin{Lemma}\label{irredcubic} In the above setting, if $\delta=3$ then $Z_s$
cannot be an irreducible cubic.\end{Lemma}

\begin{proof}  For a general point point $y\in Z_s$ we
can consider the plane $\alpha_y$, which, by Lemma \ref {tridimensional} is
different from $\alpha=\alpha_x$. \par

Notice that there is a unique point $w_y\neq y$ belonging to $T_yZ_s\cap Z_s$. Since
by Lemma \ref{Z1} one has $T_yX\cap X\subset \alpha_y$ and since both $y,w_y$ are
points of $T_yX\cap X$, we have that $T_yZ_s\subset \alpha_y$. Since one has also
$T_yZ_s\subset \alpha$ and $\alpha_y\not=\alpha$,
then $T_yZ_s=\alpha_y\cap\alpha$. \par

We also remark that $w_y\in Z_y$ and that $\alpha_y\not=\alpha_{w_y}$. Indeed,
$w_y$ is also a general point of $Z_s$ and therefore, as we saw,
$T_{w_y}Z_s=\alpha_{w_y}\cap\alpha$. Since $T_{w_y}Z_s\not=T_yZ_s$, the assertion
follows.\par

Consider now the
Zariski closure  ${\cal F}_x$ in ${\cal F}$ of the planes $\alpha_y$ with $y\in
Z_s$ a general point. There are two possibilities:
\begin{itemize}
\item[i)] the general hyperplane
through $T_xX$ contains some general plane in ${\cal F}_x$
\item[ii)] $T_xX$ and all the
planes in ${\cal F}_x$ sit in one
and the same $\Proj^4$ which we will denote by $F_x$.
\end{itemize}

Assume we are in case $(i)$ and let us keep the notation introduced at the
beginning of this section. The surface $\Sigma$ contains the curve $Z_y$ for $y\in
Z_s$ a sufficiently general point. Thus the support $Z'$ of $Z_y$ is still an
irreducible cubic. We let ${\cal Z}'$ be the strict transform of $Z'$ on $S$. It is
clear that $Z'$ is contracted by the morphism $\phi_{\cal L}$, thus its arithmetic
genus should be $0$. On the other hand, the plane cubic $Z'$ does not contain $x$,
hence ${\cal Z}'$ is isomorphic to $Z'$ and therefore its arithmetic genus is $1$, a
contradiction.\par

Assume we are in case $(ii)$. First notice that $F_x \cap X$ has dimension
$2$. Indeed it contains the irreducible surface $E_x$ described by the cubic
curves $Z_y$ with $y$ a variable point in $Z_s$. \par

Consider the Zariski
closure $V$ of the union of $F_y$ as $y\in Z_s$ varies. Remark that for $y$ general
in $Z_s$ one has $F_y=<T_yX,\alpha_{w_y}>$ hence $<\alpha_y,\alpha_{w_y}>\subset
F_y\cap F_x$, thus $\dim (F_x\cap F_y)\geq 3$.  Therefore, if we consider the
projection $\f:\Proj^7\flip \Proj^3$ from $T_xX$, then
$\dim(\f(V))=2$. On the other hand $\f(X)=\Proj^3$, hence $V$ does not contain
$X$, and this implies that the surface $E:=E_y$ does not depend on $y\in
Z_s$. This means that either $E$ has a $2$-dimensional family of plane
cubics or that the planes $\alpha_y$,  with $y\in Z_s$, pairwise meet along a
variable line. In either case $E$ has to be a cubic surface in a $\Proj^3$. By
Proposition \ref {Prop:hyperquadrics} we have a contradiction again.\end{proof}

Now we go back to the characteristic system  ${\cal C}={\cal Z}+{\cal
L}$ of $\Lambda$ on $S$. By the previous
lemmas, any irreducible component $\cal Z'$ of $\cal Z$ is
rational and smooth.  We have the
following elementary fact.

\begin{Lemma} \label{Z2} If $\cal Z'$ is a reduced, connected subcurve
of $\cal Z$, then $({\cal
L}+{\cal Z'})\cdot {\cal Z'}\leq -1+p_a({\cal Z'})$ and ${\cal Z'}^2\leq
-1+p_a({\cal
Z'})$. \end {Lemma}

\begin {proof} Consider the exact sequence

$$0\to {\cal O}_S({\cal L})\to {\cal O}_S({\cal L}+{\cal Z'})\to {\cal
O}_{{\cal Z'}} ({\cal L}+{\cal Z'})\to 0$$

\noindent Since $h^0(S,{\cal O}_S({\cal L}))=h^0(S,{\cal O}_S({\cal
L}+{\cal Z'}))$ and $h^1(S,{\cal O}_S({\cal
L}))=0$, the sequence yields $h^0({\cal Z'},{\cal O}_{{\cal Z'}} ({\cal L}+{\cal
Z'}))=0$. Hence $\deg({\cal L}+{\cal
Z'}_{|{\cal Z'}})\le -1+p_a({\cal Z'})$. The proof of the second
assertion is similar. \end{proof}

We need two more results that will be proved later on.

\begin{Proposition}\label{segrecone} Let $\tau: \Sigma\map {\Proj^{r-4}}$ be the projection of
$\Sigma$, the general hyperplane tangent section,  from $T_xX$. Then $\tau$
contracts $L$ to a point.\end{Proposition}

Finally a technical proposition on the normal sheaf of a reduced curve in $X$, for
more detail see the next section.

\begin{Proposition}\label{normal} Let $C\subset X$ be any reduced
curve with {\it superficial singular points} only, i.e. points of
embedding dimension $1$ or $2$. \par
Let $f: \tilde C\to C$ be the normalization morphism. We set
$p_n(C)=p_a(\tilde C)$. Then we have a natural exact sequence
$$0\to T\to N_f\to N'_f\to 0,$$
where $N'_f$ is locally free of rank $n-1$. We denote by $t_C$
the degree of $T$, which can be called the
{\it number of cusps} of $C$. Then
$$c_1(N_{C|X})=c_1(N'_f)+t_C+2(p_a(C)-p_n(C)).$$
\end{Proposition}

We can now prove the key Lemma in a strengthened form.

\begin{Lemma}\label{lines} If $\dim(Z)=1$,
then $Z_s$ is a bunch of $k$ lines through $x$. For each irreducible
component ${\cal Z}'$ of $\cal Z$ one has $K_S\cdot {\cal Z}'=-1$, and $K_S\cdot
{\cal Z}=-k$.\end{Lemma}

\begin{proof} We separately discuss the various possibilities for $Z_s$
given by the previous lemmas.

\begin {case}\label {case1} $Z_s$ is a bunch of lines through $x$.\end
{case}

Let $Z_i$ be the lines in $Z_s$ and ${\cal Z}_i$ their strict transform
on $S$, with $i=1,...,k$.
We have ${\cal Z}=a_1{\cal Z}_1+...+a_k{\cal Z}_k$, with $a_i,
i=1,...,k,$ positive integers. By
Proposition \ref {segrecone} we have ${\cal L}\cdot {\cal Z}_i=0, {\cal
Z}_i^2\leq -1$ for $i=1,...,k$. Hence
$1=H'_{|S}\cdot {\cal Z}_i=\sum_{j=1}^k a_j{\cal Z_i}\cdot {\cal
Z_j}+2=a_i{\cal Z}_i^2+2$, thus $a_i=1$ and
${\cal Z}_i^2=-1$, whence the assertion.

\begin {case}\label{case3} $Z_s$ is an irreducible conic.\end{case}

If ${\cal Z'}$ is the strict transform of $Z_s$ on $S$, we have ${\cal
Z}=a{\cal Z'}$, with $a$ a positive number.
 Arguing as before, we have $2={\cal Z'}\cdot ({\cal L}+a{\cal
Z'}+2F)\leq 1+ (a-1){\cal Z'}^2$, which
leads to a contradiction by Lemma \ref{Z2}.

\begin {case}\label{case4} $Z_s$ is a reducible conic.\end{case}

Then $Z_s$ consists of two lines $Z_1, Z_2$ which we may assume not both
passing through $x$. We denote by ${\cal
Z}_1, {\cal Z}_2$ the strict transforms of $Z_1, Z_2$ on $S$. Notice
that ${\cal Z}_1\cdot {\cal
Z}_2=1$. Moreover ${\cal Z}=a_1{\cal Z}_1+a_2{\cal Z}_2$, with $a_1,a_2$
positive integers. \par

Set ${\cal Z}'={\cal Z}_1+{\cal Z}_2$. By Lemma \ref{Z2} we have
$({\cal L}+{\cal Z}')\cdot {\cal Z}'\leq -1$.
Hence $2\leq (a_1-1){\cal Z}_1^2+(a_2-1){\cal Z}_2^2+(a_1+a_2-1)$, a
contradiction.

\begin {case}\label{case5} $Z_s$ is a plane cubic which is the union of
an irreducible conic and of a line.\end
{case}

Let $Z_1$ be the line and $Z_2$ the conic in $Z_s$ and let ${\cal Z}_1$,
${\cal Z}_2$ be the strict transforms
on $S$ of $Z_1$ and $Z_2$ respectively. Notice that, by Lemma \ref
{le:kfil}, $x$ belongs to $Z_1$. Furthermore $Z_1$ is at least $2$-filling
and, by Lemma \ref {tridimensional}, $Z_2$ is at least $3$-filling for $X$.\par

Suppose both $Z_1$ and $Z_2$ contain $x$. Then $Z_s$ is a Cartier
divisor
on $\Sigma$, and one has $c_1(N_{Z_s|X})=Z_s^2+3$. Consider the
normalization map
$f: \Proj^1\amalg\Proj^1\to Z_s$. By Lemma \ref{filling} we have
$c_1(N'_f)\geq  1$, hence by Proposition \ref
{normal} we have $c_1(N_{Z_s|X})\geq 5$. Accordingly one has
$Z_s^2\geq 2$ and therefore $({\cal Z}^{\prime})^2\geq 0$, a contradiction.\par

Suppose $Z_1$ contains $x$ and $Z_2$ does not. Then $Z_2$ is a Cartier
divisor on $\Sigma$ and, arguing as before
we see that ${\cal Z}_2^2=-1$, thus ${\cal L}\cdot {\cal Z}_2=0$ by
Lemma \ref {Z2}. Furthermore ${\cal L}\cdot
{\cal Z}_1=0$ by Proposition \ref {segrecone}. Since ${\cal Z}_1\cdot
{\cal Z}_2=2$ we have a contradiction
to the last assertion of Proposition \ref {Le:br}.\par

\begin {case}\label {case6} $Z_s$ is a plane cubic which is the union of
three lines, not all passing through
$x$.\end {case}

Let $Z_i$ be these lines, and let ${\cal Z}_i$ be the strict transform
of $Z_i$ on
$S$, $i=1,2,3$. By Lemma \ref {le:kfil}, $x$ has to lie on at least
one of the lines $Z_i$, $i=1,2,3$.

Suppose that $x\in Z_1$ but $x\notin C:=Z_2\cup Z_3$. Then $C$ is a
Cartier divisor on $\Sigma$ and one has
$c_1(N_{C|X})=C^2+2$. Consider the normalization map $f:
\Proj^1\amalg\Proj^1\to C$. One has $c_1(N'_f)\geq  0$,
hence by Proposition \ref {normal} we have $c_1(N_{C|X})\geq 2$, thus
$C^2\geq 0$ which leads to a contradiction to
Lemma \ref {Z2}.\par

Suppose that $x\notin Z_1$ but $x\in Z_2\cap Z_3$. Again we set $C:=Z_2\cup Z_3$.
As before, we see that $C^2\geq 0$. Hence $({\cal
Z}_2+{\cal Z}_3)^2\geq -2$. Thus ${\cal Z}_2^2={\cal Z}_3^2=-1$.
Moreover, by a
similar argument, we have that ${\cal Z}_1^2=-1$. Thus ${\cal L}\cdot
{\cal Z}_i=0, i=1,2,3$. \par

Let ${\cal C}={\cal L}+a_1{\cal Z}_1+a_2{\cal Z}_2+a_3{\cal Z}_3$, with
$a_i$ positive integers, $i=1,2,3$. By intersecting $H_{|S}'$ with $Z_i,
i=1,2,3$, we find the relations
$a_1=a_2+a_3-1=a_2-1=a_3-1$, which lead to $a_2=a_3=0, a_1=-1$, a
contradiction. \end{proof}\medskip

We are finally in a position to prove Theorem \ref{gradoOADP}

\begin{proof}[Proof of Theorem \ref{gradoOADP}]

As already noticed, we may assume that $X$ is not a scroll over a curve. By
Riemann-Roch theorem, one has:

\begin{equation} \label{circ} 3-h^1(S, {\cal O}_S({\cal C}))=\chi({\cal O}_S({\cal
C}))=1+\frac12({\cal C}^2-K_S \cdot {\cal C}).\end{equation}

\noindent One has ${\cal C}^2=\deg(X)-8$ and by Lemma \ref{Le:br} one
has
$K_{S}\cdot \L=-3$. Let $k$ be the number of lines through $x$. By Lemma \ref{lines}
we have $K_S\cdot
{\cal C}= -3-k$. Finally, by
equation (\ref{circ}), we have $\deg(X)=9-k$.\end{proof}

\begin{remark} \label{rette} The proof of Theorem \ref{gradoOADP} shows
that actually $d\leq 8$, unless for a
general point $x\in X$, the scheme $Z_x$  cut out by $T_xX$ on $X$ is
$0$-dimensional supported at $x$ (see Proposition \ref {Le:br}). Furthermore, if
$X$ is not a scroll, i.e. $d\geq 6$, and if $d\leq 8$ then $Z$ consists of $9-d$
lines through $x$. In particular in these
cases $X$ is ruled by lines. We remark that this is in fact the
situation in the known examples of Edge varieties,
$d=6,7$ and of the scroll of degree $8$ (see examples \ref{edgevarieties}  and
\ref {degree8}).\end{remark}

\section {A few technical results}\label{lemmata}

In this section we prove a few technical results we have been using in
the previous section.\medskip

First we  sketch
the proof of Proposition \ref{scroll} for lack of a modern reference.

\begin{proof}[Proof of Proposition \ref{scroll}] We prove only the non trivial implication. Assume that
$n<r-1$. By cutting $X$ with
a general hyperplane we can reduce ourselves to the surface case $n=2$.
The general tangent
hyperplane section $H$ of $X$ is reducible. If $H$ is non-reduced, then
the dual variety of $X$ is degenerate and
$X$ is a developable scroll (see \cite {GH}). Otherwise write $H=A+B$.
Since $H$ is a nodal curve, with a single
node at the point $x$ of tangency, and since $H$ is connected, we see
that $A$ and $B$ are smooth curves meeting
at $x$. Namely $H$ is $1$-connected but not $2$-connected. This implies
that either $A$ or $B$ is a line and $X$
is a scroll (see \cite {vdv}). \end {proof}

Notice that, in the case $X$ is smooth, as remarked in the proof of
Proposition \ref {RPT}, the previous result
follows from \cite {Kl} and \cite[Theorems 2.1 and 3.2]{Ei}. As a consequence we have
the following.

\begin{Proposition} \label{OADPrigate} Let $X\subset\Proj^{2n+1}$ be a
{\it OADP} $n$-fold. Let $x$ be a
general point of $X$. The
intersection scheme of $X$ with $T_xX$ has an irreducible component of
dimension $n-1$ if and only if $X$ is one of the two scrolls
$S(1^{n-1},3)$ and $S(1^{n-2},2^2)$.\end{Proposition}

\begin{proof} This is a consequence of Lemma \ref {scroll} and
Proposition \ref {OADPscrolls}. \end
{proof}\medskip

Next we justify the basic proposition about linear systems of rational curves on surfaces.

\begin{proof}[Proof of Proposition \ref{rational}]
Assume $|\cal M|$ has some base point. There is no lack
of generality in supposing that it
has only one base point at $x\in S$ which is an ordinary multiple point
of multiplicity $m$ for the general curve
$C$ in $|\cal M|$. Indeed, we can always put ourselves in this situation
after having performed a suitable
sequence of blow-ups. \par

Let $f: S'\to S$ be the blow-up of
$S$ at $x$ with exceptional divisor $E$, and consider the proper
transform $C'$ of $C$ on $S'$, i.e.
$C'=f^*(C)-mE$. Note that $C'$ is a smooth curve of genus $g$ on $S'$.
If we set $f^* ({\cal M}):={\cal M}'$, the
curve $C'$ determines the line bundle ${\cal M}'\otimes {\cal
O}_{S'}(-mE)$. Look at the exact sequence:

\begin{equation}\label{plus} 0\to {\cal O}_{S'}\to {\cal O}_{S'}(C')\to {\cal
O}_{C'}(C')\to 0\end{equation}

\noindent Since $h^1(S',{\cal O}_{S'})=0$, one has $h^0(C',{\cal
O}_{C'}(C'))=h^0(S',{\cal
O}_{S'}(C'))-1=h^0(S,{\cal M})-1$. Thus $h^0(C',{\cal O}_{C'}(C'))\geq
g+1$ and therefore $h^1(C',{\cal
O}_{C'}(C'))=0$. Hence $g+1\leq h^0(C',{\cal O}_{C'}(C'))=(C^{\prime})^2-g+1$,
i.e.
$(C^{\prime})^2\geq 2g$ and then $K_{S'}\cdot C'<0$,  proving that all the
plurigenera of $S'$ vanish. Since the same is
true for $S$, we have that $S$ is rational.  Now, from equation (\ref{plus}) we deduce
$h^1(S', {\cal M}'\otimes {\cal O}_{S'}(-mE))=h^1(S',{\cal
O}_{S'}(C'))=0$. Look at the exact sequence:

$$0\to {\cal M}'\otimes {\cal O}_{S'}(-mE) \to {\cal M}'\to {\cal
M}'\otimes {\cal O}_{mE}\simeq {\cal
O}_{mE}\to 0$$

\noindent Since $h^0(S', {\cal M}'\otimes {\cal O}_{S'}(-mE))=h^0(S',
{\cal M}')$, we see that $h^1(S',
{\cal M}'\otimes {\cal O}_{S'}(-mE)))\geq h^0(S',{\cal M}'\otimes {\cal
O}_{mE})=m(m+1)/2>0$, a
contradiction. \end{proof}\medskip

Here we introduce in a more detailed way the notion of {\em number of cusps} for a curve $C$ with
only superficial singular points.

Let now  $X$ be a smooth, projective variety of dimension
$n$ and let $C\subset X$ be a reduced
curve with {\it superficial singular points} only, i.e. points of
embedding dimension $1$ or $2$. \par

First we remark that the normal sheaf $N_{C|X}$ of $C$ in $X$ is locally
free of rank $n-1$ (see \cite[Lemma 2.1]{MP}). \par

Let $f: \tilde C\to C$ be the normalization morphism. We set
$p_n(C)=p_a(\tilde C)$. We can also consider the
normal sheaf to the map $f$, denoted by $N_f$, which is defined by the
exact sequence:
$$0\to T_{\tilde C} \buildrel df \over\rightarrow f^*(T_X)\to N_f\to 0$$
 The sheaf $N_f$ is locally free, except at
the points of $\tilde C$ where the differential of $f$ vanishes, where
$N_f$ has torsion. If $T$ is the torsion
subsheaf of $N_f$, we have the sequence
$$0\to T\to N_f\to N'_f\to 0$$
 where $N'_f$ is locally free of rank $n-1$. We denote by $t_C$
the degree of $T$, which can be called the
{\it number of cusps} of $C$. \par

\begin{proof}[Proof of Proposition \ref{normal}]
Notice that there is a smooth surface $S\subseteq X$
containing $C$. We denote by $\phi$
the map $\phi: \tilde C\to S$ induced by $f$. We have the sequences:

$$0\to N_{C|S}\to N_{C|X}\to N_{{S|X}_{|C}}\to 0$$

$$0\to N'_\phi\to N'_f\to f^*N_{{S|X}_{|C}}\to 0$$

\noindent which shows that we can reduce ourselves to the case $X$
itself is a surface. Then, by definition, one
has:

$$c_1(N'_f)=-\deg(f^*(K_X))+2p_a(\tilde C)-2-t_C$$

\noindent and by adjunction:

$$c_1(N_{C|X})=-K_X\cdot C+2p_a(C)-2$$

\noindent which proves the assertion.\end{proof}\medskip

%\subsection{}\label{se:filling}
We are left to prove Proposition \ref{tangentplanes} and Proposition
\ref{segrecone}. These we believe  to be of
independent interest.

Let $C$ be an irreducible reduced rational curve, so that
$\tilde
C\simeq \Proj^1$. Then $N'_f\simeq {\cal
O}_{\Proj^1}(a_1)\oplus...\oplus {\cal O}_{\Proj^1}(a_{n-1})$, with
$a_1\leq ...\leq a_{n-1}$.

Under this hypothesis the above construction gives.

\begin{Lemma}\label{filling} If $C$ is a rational $k$-filling curve for
$X$ then $a_1\geq 0$ and
$a_1+...+a_{n-1}\geq k-n+1$. In particular $h^1(\tilde{C},N_f)=h^1(\tilde{C},N'_f)=0$.  \end{Lemma}

\begin{proof} A count of parameters shows that $k\geq n-1$. Hence we
have a dominant morphism ${\Proj^1}\times B\to
X$ where $B$ is a suitable $(n-1)$-dimensional variety. In this
situation, given a general point $b\in B$, there
is a {\it focal} square matrix $\Phi$ arising (see \cite {cs}), whose
rows
are given by the sections of $N'_f$ corresponding to the $n-1$
infinitesimal deformations of $f: {\Proj^1}\times
\{b\}\to X$ determined by $n-1$ independent vectors of $T_b(B)$. By the
filling hypothesis, $\Phi$ has maximal
rank, which implies $a_1\geq 0$. The final assertion translates the fact
that $k\leq h^0(\tilde C,N_f')$ (see \cite
{AC}).\end{proof}\medskip

%\subsection

Then a lemma about foci of certain families of planes in $\p^4$ (again,
see \cite{cs} for the basics in the theory of foci).

\begin{Lemma}\label{tangfocal} Let $X\subset \Proj^4$ be a non degenerate, irreducible
hypersurface. Suppose there is an irreducible subvariety $\cal F\subset
\G(2,4)\times \p^4$ of dimension $2$, such that the general point $\xi\in\cal F$
corresponds to a pair $(\Pi,p)$ where $\Pi$ is a plane in $\Proj^4$, $p$ is
a smooth point of $X$ and $\Pi$ is tangent to $X$ at $p$. Suppose the
projection $\cal G$ of $\cal F$ to $\G(2,4)$ has dimension $2$. Then $p$ is a focus
for the family $\cal G$ on $\Pi$.\end{Lemma}

\begin{proof} Let ${\bf x}=(x_0,...,x_4)$ be homogeneous coordinates on $\p^4$, let
${\bf a}=(a_0,...,a_4)$ be the dual coordinates and let $(u_1,u_2)$ be local
parameters on $\cal F$ around $\xi$. Thus we may assume that $(u_1,u_2)$ are local
parameters on $\cal G$ around $\Pi$, so that the equations of $\Pi$ are given by:

\begin{equation} \label{pplanes}  {\bf a}_i\times {\bf x}=0, i=1,2 \end{equation}

\noindent where ${\bf a}_i={\bf a}_i(u_1,u_2)$ are regular functions of $(u_1,u_2)$.
We will denote by a subscript $j$ the differentiation with respect to $u_j$,
$j=1,2$. Notice that the foci of $\cal G$ on $\Pi$ are defined by the additional
equation:

\begin{equation} \label{pfoci} \det({\bf a}_{i,j}\cdot {\bf x})_{\{i,j=1,2\}}=0
\end{equation}

We can now assume that the homogeneous coordinates ${\bf p}=[p_0,...,p_4]$ of $p$ are
also regular functions ${\bf p}(u_1,u_2)$ of $(u_1,u_2)$. Hence, by differentiating
the equations (\ref{pplanes}) with ${\bf x}={\bf p}$, we get:

\begin{equation} \label{relations}
{\bf a}_{i,j}\times {\bf p}+{\bf a}_i\times {\bf p}_j=0, i,j=1,2  \end{equation}

Let

\begin{equation} \label{eq} f({\bf x})=0\end{equation}

\noindent  be the equation of $X$. The
equation of $T_pX$ is given by:

\begin{equation}\label{tangent} \sum_{h=0}^4 {\frac{\partial f}{\partial
x_h}}({\bf p})x_h=0\end{equation}

\noindent By differentiating equation (\ref{eq}) with ${\bf x}={\bf p}$, we also get:

\begin{equation}\label{der} \sum_{h=0}^4{\frac{\partial f}{\partial x_h}}({\bf
p})p_{h,j}=0, j=1,2\end{equation}

Let ${\bf t}=[t_1,t_2]$ be homogeneous coordinates on $\p^1$. By the hypothesis,
$T_pX$ contains $\Pi$. Thus there is a regular function ${\bf t}(u_1,u_2)$ of
$(u_1,u_2)$ such that the equation (\ref{tangent}) of $T_xX$ is also given by:

$$(t_1{\bf a}_1+t_2{\bf a}_2)\times {\bf x}=0$$

\noindent By taking into account equation (\ref{der}), we find:

$$(t_1{\bf a}_1+t_2{\bf a}_2)\times {\bf p}_j=0, j=1,2$$

\noindent  and using equation (\ref{relations}) we get:

$$(t_1{\bf a}_{1,j}+t_2{\bf a}_{2,j})\times {\bf p}=0, j=1,2$$

\noindent The assertion follows by the equation \ref{pfoci} for the focal locus on
$\Pi$.\end{proof}

Now we can prove our first proposition.

\begin{proof}[Proof of Proposition \ref{tangentplanes}] Suppose there is an irreducible
subvariety $\cal F\subset \G(2,r)$ of dimension $2$,  such that the
general point $p\in \cal F$ corresponds to a
plane $\Pi\subset \Proj^r$ which is tangent to $X$ along a curve
$\Gamma$ which is $2$-filling. Make a general
projection in $\Proj^4$. We thus get a hypersurface $X'$ in $\Proj^4$,
the projection of $X$, such that there is an
irreducible subvariety $\cal G\subset \G(2,4)$ of dimension $2$,  such
that the general point $p'\in \cal G$
corresponds to a plane $\Pi'\subset \Proj^4$, the projection of $\Pi$,
which is tangent to $X'$ along a curve
$\Gamma'$, the projection of $\Gamma$.  By Lemma \ref{tangfocal}, every
point of $\Gamma'$ is a focus for the family $\cal
G$ on $\Pi'$. Hence
$\deg(\Gamma)=\deg(\Gamma')\leq 2$. \par

Assume that $\deg(\Gamma)=1$. By lemma
\ref {filling}, we have $N_{\Gamma|X}={\cal
O}_\Gamma(a)\oplus {\cal O}_\Gamma(b)$ with $a,b\geq 0$ and $\Gamma$ is
unobstructed.  By the results of \cite {dg}, we have an injection ${\cal
O}_\Gamma(1)\to
N_{\Gamma|X}$, which implies $a+b\geq 1$. Then $X$ has
a $3$-dimensional family of lines at least, hence it is either
$\Proj^3$, or a quadric
in $\Proj^4$ or a scroll
over a curve (see \cite  {rogora}), a contradiction.\par

Assume now $\deg(\Gamma)=2$. We may suppose $\Gamma$ to be irreducible,
hence $\Gamma\simeq \Proj^1$. Set $N_{\Gamma|X}={\cal
O}_{\Proj^1}(a)\oplus {\cal O}_{\Proj^1}(b)$. As above, we have $a,b\geq
0$ and $\Gamma$ is unobstructed. If
$a+b\geq 4$ then $X$ has a $6$-dimensional family of conics at least.
This means that its general hyperplane section
has a $3$-dimensional family of conics (see \cite[Lemma 2.3]{MP}),
hence $X$ is a quadric in $\Proj^4$, a
contradiction. By the results of \cite {dg} we have an injection ${\cal
O}_\Gamma(2)\to N_{\Gamma|X}$, which
is incompatible with $a+b\leq 3$.\end{proof}

\begin{remark} If $X$ is a smooth quadric in $\Proj^4$, then for every
line $L$ on $X$ there is a unique plane
$\Pi$ tangent to $X$ along $L$. Namely $\Pi=\cap_{x\in L}T_xX$.\par

On the other hand one may prove that the hypothesis $X$ not a scroll
over a curve, in the previous
proposition, is unnecessary. We will not dwell on this here.
\end{remark}
\medskip

Let
$X\subset {\Proj^r}$, $r\geq 7$, be a smooth,
irreducible projective $3$-fold which is not a scroll over a curve. Suppose that $X$
is ruled by lines. Let $x\in X$ be a general point and let $L$
be a line through $x$. Let $\Sigma$ be a general tangent hyperplane
section of $X$ at $x$. Notice that $\Sigma$
contains $L$. Let $\tau: \Sigma\map {\Proj^{r-4}}$ be the projection of
$\Sigma$ from $T_xX$. Note that $\Sigma$
has a single node at $x$ (see \cite[Theorems 2.1]{Ei}), thus it is
normal and therefore $\tau$ is well
defined at the general point of $L$.

We pay all our debits with the following enlarged version of Proposition \ref{segrecone}.
\begin{Proposition} In the above setting, $\tau$
contracts $L$ to a point.\end{Proposition}

\begin{proof}
 Let $y\in L$ be the general point of $L$. It is clear that
$\tau(y)$ is the
intersection of the target ${\Proj^{r-5}}$ with the $\Proj^4$ spanned by
$T_xX$ and by $T_y\Sigma$. If $H$ is the
hyperplane cutting $\Sigma$ on $X$, one has $T_y\Sigma=T_yX\cap H$. \par

The variety $S_{L,X}=\cup_{y\in L}T_yX$ is the so-called {\it Segre
cone} of $X$ along $L$ (see \cite[chapter 1, \S 3.1]{rogora}).
This is a quadric cone with vertex $L$ lying in a
$\Proj^5$ containing $L$ (see \cite[chapter 1, \S 3.4]{rogora}). A general tangent hyperplane $H$ at $x$
cuts $S_{L,X}$ along the union of $T_xX$
with another $3$-space $\Pi$, intersecting $T_xX$ along a plane
containing $L$. Thus $\Pi$ and $T_xX$ span a
$\Proj^4$, whose projection from $T_xX$ is a fixed point. This proves
that $\tau(y)$ stays fixed as $y$ varies on
$L$.\end{proof}

\section{Classification of smooth 3-folds with one apparent double
point}\label {classification}

In this section we  apply the previous results to classify smooth
{\it OADP} 3-folds. Our main result is the following theorem:

\begin{Theorem} Let $X\subset\Proj^7$ be an {\it OADP} 3-fold.  Let $d$
be the degree of $X$.
Then:
\begin{itemize}
\item[i)] either $d=5$ and $X$ is a scroll in planes of type
$S(1^2,3)$ or $S(1,2^2)$ as in Example \ref{scrolls},
\item[ii)] or $d=6,7$ and $X$ is an Edge variety as in
Example \ref{edgevarieties},
\item[iii)] or $d=8$ and $X$ is a scroll in lines over $\Proj^2$ as
described in
Example \ref {degree8}.\label{Th:classification}
\end{itemize}
\end{Theorem}

Since we have shown before that {\it OADP} 3-folds
have degree $d$ less than or equal to 9, then these varieties
are contained in the lists of \cite{Io1}, \cite{Io2} and
\cite{Io4}. Hence it could be possible to deduce Theorem \ref
{Th:classification} from Ionescu's results.
However we have already described, on the way,  several
geometrical properties
of {\it OADP} 3-folds. This simplifies a lot the classification. To finish the proof
of the classification we only need a few basic properties of generalized adjunction
maps and some of the ideas contained in \cite{Io1}. \par

As we saw, $X$ is linearly normal and regular. Then its general curve
section is a linearly normal curve $C\subset \Proj^5$ of degree $d$
and sectional genus $g$.  By Clifford theorem, the curve $C$ is non
special, hence
by Riemann-Roch theorem
we have $d=g+5$ and hence $g\leq 4$.\par

>From now on we can suppose $6\leq d\leq 9$, or equivalently
$1\leq g\leq 4$ and $X$ not a scroll over a curve.
Moreover, by \cite[Theorem 1.4]{Io1}, we can suppose that
$|K_X+2H|$ is base point free and that it gives a morphism
$\phi=\phi_{|K_X+2H|}:X\to\Proj^{g-1}$, the so-called {\it adjunction
map} of the polarized pair
$(X,H)$.

The following proposition will conclude the classification.

\begin{Proposition} Let $X\subset\Proj^7$ be a smooth {\it OADP} 3-fold
of degree $d\geq 8$. Then $d=8$ and $X$ is like in
Example \ref {degree8}.\label{Le:89}
\end{Proposition}

\begin{proof} Assume $d=8$ and therefore $g=3$. Let $\phi:X\to\Proj^2$
be the adjunction
map. Since $\phi(X)$ is non-degenerate, then either $\phi(X)$ is a
curve or $\phi(X)=\Proj^2$. If $\phi(X)$ were a curve, then
$X$ would be a hyperquadric fibration by
\cite[Proposition 1.11]{Io1} and  by Proposition \ref{Prop:hyperquadrics} this is
impossible. Thus $\phi:X\to\Proj^2$ is surjective and endows $X$ with a
structure of ${\Proj^1}$-bundle over
${\Proj^2}$ (see again \cite {Io1}). \par

 We  now  show that $X$ is like in Example \ref{degree8}. Let $H$ be
 a general hyperplane section of $X$. Then applying $\phi_*$ to
the exact sequence

$$0\to\O_X\to\O_X(H)\to\O_H(H)\to 0$$

\noindent we
obtain

$$0\to\O_{\Proj^2}\to {\E}\to\phi_*(\o_H(H))\to
R^1\phi_*(\O_X)=0,$$

\noindent where $\E$ is locally free of rank 2 and $X=\Proj({\cal
E})$.\par

Since $C=H^2$ is a smooth curve of genus 3  and degree 8 on
$H$, which is also smooth, we have that
$C=\phi^*(\O_{\Proj^2}(4)-E_1-\ldots-E_8)$, where
$p_i=\phi(E_i)$ are 8 distinct points. Of course $E_i$ is a $(-1)$-curve
and
$C\cdot E_i=1$, for all $i=1,...,8$. It follows that
$\phi_*(\o_H(H))=\I_{p_1,\ldots, p_8}(4)$. Since $H$ is smooth, the 8
points do
not lie on a conic and no four point are collinear. Moreover, we  claim
that, by the
generality of $H$, we can suppose that no seven points $p_i$ lie on a
conic. Once this claim is proved we can conclude that $X$ is like in
Example \ref{degree8}. \par

To prove the claim we recall that $Pic(X)=<H,
\phi^*(\o_{\Proj^2}(1))>$. If the claim were not true, for the general
$H$
we would have a line $l_H\subset H$, the strict transform on $H$
of the conic through the seven points $p_i$. Notice that
there would be a $r$-dimensional family of lines $l_H$ on $X$, with
$r\geq 2$.
Clearly
$H\cdot l_H=1$ and $\phi^*(\o_{\Proj^2}(1))\cdot l_H=2$. Remark that
$h^0(\o_X(H)\otimes\phi^*(\o_{\Proj^2}(-1)))=h^0(\E(-1))=2$, thus the
system $|H+\phi^*(\o_{\Proj^2}(-1))|$
is a pencil. Take a divisor
$D\in|H+\phi^*(\o_{\Proj^2}(-1))|$.  From $D\cdot l_H=-1$ we conclude
that the lines $l_H$ cannot fill up $X$,
hence they fill up a plane $\Pi$ which sits in the base locus $F$ of the
pencil $|H+\phi^*(\o_{\Proj^2}(-1))|$
and surjects onto ${\Proj^2}$ via $\phi$. Call $|M|$ the movable part of
the pencil.  Set
$l_p=\phi^{-1}(p)$ for a $p\in \Proj^2$. The relations $1=D\cdot
l_p=M\cdot l_p+F\cdot
l_p$ and  $F\cdot l_p>0$  would imply $M\cdot l_p=0$ so that
$M=\phi^*(\o_{\Proj^2}(a))$, $a>0$. Then we would have $2=h^0(X,{\cal
O}_X(M))=h^0({\Proj^2}, \o_{\Proj^2}(a))$,
which is impossible.\par

Let us now suppose $d\geq 9$. Then by Theorem \ref{gradoOADP} $d=9$ and therefore $g=4$. Then $1\leq
\dim(\phi(X))\leq 3$.\par

If $\dim(\phi(X))=1$, then $X$ would be a
hyperquadric fibration and this is in contrast with
Proposition \ref{Prop:hyperquadrics}.\par

 If $\dim(\phi(X))=2$, then
the adjunction map would determine a structure of scroll in lines over a smooth
surface, see {\it loc. cit.}. This is not possible by Theorem \ref {gradoOADP}.\par

Let us now exclude the case $\phi(X)=\Proj^3$. Since the sectional genus is $g=4$
in this case, the only possibility left out by
\cite[Theorem II]{Io3} is that $X$ could be the blow-up at finitely
many points of a  $\Proj^2$-bundle  over $\Proj^1$, which we denote by $X'$.
Moreover $\phi$ would factor through $X'$ and, by
\cite[Lemma 1.2]{Io3} we would have a finite
morphism $\phi':X'\to \Proj^3$ embedding   each  fiber of the projection $X'\to
\Proj^1$ as a plane, which is
clearly impossible.
\end{proof}
\medskip

We are finally ready for the proof of our classification
theorem.\medskip

\begin{proof}[Proof of Theorem \ref{Th:classification}]  We know by
theorem
\ref {gradoOADP} that $5\leq d\leq 9$. If $d=5$, then $X$
is either $S(1,2^2)$ or $S(1^2,3)$ by lemmas \ref{basiscrolls} and
\ref{OADPscrolls}.

Then we can suppose $6\leq d\leq 9$, or equivalently
$1\leq g\leq 4$ and $h^0(K_X+2H)=g$. Consider the adjunction morphism
$\phi=\phi_{|K_X+2H|}:X\to\Proj^{g-1}$.

If $g=1$, then $K_X=-2H$ and $X$ is a Del Pezzo manifold of
degree 6. Iskovskikh has shown that either $X$ is the Segre variety
$Seg(1^3)$, i.e. the Edge 3-fold of degree 6
(see Example \ref {edgevarieties}), or $X=\p(T_{\Proj^2})$, see
\cite{Is1} and \cite{Is2}. The last variety is the
hyperplane section of $Seg(2^2)$ and  its secant variety is  a
hypersurface of degree 3,
i.e. $\p(T_{\Proj^2})$ is defective and therefore is not an {\it
OADP}-variety (see \cite {Zak1}, \cite {cc}).

If $g=2$, then  $\phi:X\to\Proj^1$ gives  to $X$ the structure of
hyperquadric fibration over $\Proj^1$, see \cite[Proposition
1.11]{Io1}. By Proposition \ref{Prop:hyperquadrics}, $X$ is  an Edge variety.

If $3\le g\leq 4$ we can apply Lemma \ref{Le:89} to
conclude.
\end{proof}

\section{Mukai varieties with one apparent double point}\label {mukai}

In this section we apply Theorem \ref{Th:classification} to
the classification of smooth $n$-dimensional {\it OADP}-varieties
$X\subset\Proj^{2n+1}$ having as a general
curve section a canonical curve $C\subset\Proj^{n+2}$. Since $X$ is
regular and linearly normal, then $C$ is linearly normal with genus
$g=n+3$. By Theorem \ref{Th:classification} we can suppose $n\geq 4$.
Moreover,  by the adjunction formula, $-K_{X}=(n-2)H$, i.e.
$X$ is a Fano variety of coindex 3. These varieties have been recently
called {\it Mukai
varieties}.  Mukai varieties have been classified in \cite{Mu}
under the assumption of the existence of a smooth divisor in
$|H|$, a condition which is clearly satisfied in our case. In any event,
this
restriction has been removed
in \cite{Me}.\par

We notice that the {\it OADP}-variety in the Example \ref {spinor} and
the one $\G^{lag}_\R(2,5)$ from Example \ref
{lagrassmann} are Mukai varieties. \par

The classification result of Mukai {\it OADP}-varieties is the
following.

\begin{Theorem} Let $X\subset\Proj^{2n+1}$ be a Mukai {\it OADP}-variety of
dimension $n$. Then either $n=4$ and
$X\subset\Proj^9$ is a linear section of the spinor  variety
$S^{(4)}\subset\Proj^{15}$ as described in Example \ref {spinor}, or
$n=6$ and
$X= \G^{lag}_\R(2,5)$.\label{Th:Mukai}
\end{Theorem}

\begin{proof} By theorems \ref {OAPDsurfaces} and
\ref{Th:classification}
we can suppose $n\geq 4$.\par

If $b_2(X)=1$, then $X$ is (a linear section of) a homogeneous
manifold as described in theorems 1 and 2 of \cite{Mu} and
necessarily $g=7,8,9$. The first and third case correspond to the cases
in the statement of the theorem. The
Mukai manifold of dimension 8 and  of sectional genus 8 is the
Grassmannian
${\mathbb{G}}(1,5)\subset \Proj^{14}$, whose secant variety is a
cubic hypersurface in $\Proj^{14}$ (see \cite {Zak1}). Therefore its
linear section with a
general $\Proj^{11}$ is  a defective variety $X_5\subset\Proj^{11}$.\par

We now show that the case $b_2(X_n)\geq 2$ is impossible. By
 \cite[Theorem 7]{Mu} we know that either $n=4$ and $X$ is a
Mukai 4-fold of product type, i.e. isomorphic to $\Proj^1\times Y$
with $Y$ a Fano 3-fold of index 2 or 4, or $X$ is (a linear
section of) one of the nine varieties described in example 2 of
\cite{Mu}. In all these cases one never has an $n$ dimensional variety
$X$ in $\Proj^{2n+1}$, except when either
$X$ is the Segre embedding of $\Proj^1\times F$, with $F$ a smooth cubic
hypersurface in $\Proj^4$, or $X$ is the
complete intersection of the cone over $Seg(2^2)$ from a point with a
quadric in $\Proj^9$. In both cases one has
$n=4$. However they do not lead to {\it OADP}-varieties. Indeed the
former variety sits in $Seg(1,4)$ as a
divisor of type $(0,3)$ and, by Proposition \ref {Le:Edge}, it is not an
{\it OADP}-variety. In the latter case $X$
is defective as well as $Seg(2^2)$ (to see this, argue as in \cite[Example 2.4]{cc},
). \end{proof}

\end{document}